\documentclass[11pt]{article}
\usepackage{amsmath}
\usepackage{dsfont}
\usepackage{amsmath}
\usepackage{amsmath,bm}
\usepackage{mathrsfs}
\usepackage{amsmath,amssymb}
\usepackage{amsfonts}
\usepackage{hyperref}
\usepackage{amsthm}
\usepackage{graphicx}
\usepackage{subfigure}
\usepackage{xcolor}
\usepackage{appendix}
\usepackage{cite}

\usepackage{tikz}

\usetikzlibrary{calc,shapes,arrows}

\usepackage{pgfplots}

\usetikzlibrary{patterns}

\usepgfplotslibrary{dateplot}

\renewcommand{\qed}{\hfill\small{$\square$}\normalsize}

\hfuzz=\maxdimen
\theoremstyle{definition}
\newtheorem{lemma}{Lemma}[section]
\newtheorem{definition}{Definition}
\newtheorem{proposition}[lemma]{Proposition}
\newtheorem{theorem}[lemma]{Theorem}
\newtheorem{corollary}[lemma]{Corollary}

\newtheorem{remark}{Remark}

\newtheorem{conjecture}{Conjecture}
\numberwithin{equation}{section}
\renewcommand{\proof}{\textbf{Proof. }}
\renewcommand{\qed}{\hfill\small{$\square$}\normalsize}

\DeclareFixedFont{\Acknowledgment}{OT1}{cmr}{bx}{n}{14pt}
\begin{document}

\title{\bf Fractional combinatorial Calabi flow on surfaces}
\author{Tianqi Wu, Xu Xu}
\maketitle

\begin{abstract}
Using the fractional discrete Laplace operator for triangle meshes,
we introduce a fractional combinatorial Calabi flow for discrete conformal structures on surfaces, which unifies and generalizes
Chow-Luo's combinatorial Ricci flow for Thurston's circle packings, Luo's combinatorial Yamabe flow
for vertex scaling and the combinatorial Calabi flow for discrete conformal structures on surfaces.
For Thurston's Euclidean and hyperbolic circle packings on triangulated surfaces,
we prove the longtime existence and global convergence
of the fractional combinatorial Calabi flow.
For vertex scalings on polyhedral surfaces, we do surgery on the fractional combinatorial Calabi flow by edge flipping
under the Delaunay condition to handle the potential singularities along the flow.
Using the discrete conformal theory established in \cite{GGLSW,GLSW},
we prove the longtime existence and global convergence of the fractional combinatorial Calabi flow with surgery.
\end{abstract}




\section{Introduction}

Since Chow-Luo's introduction  of the combinatorial Ricci flow for Thurston's circle packings on surfaces \cite{CL},
combinatorial curvature flows have been important approaches for finding geometric structures on low dimensional
manifolds, which have lots of applications in geometric topology and practical applications \cite{DGL,GY,ZG}.
The combinatorial curvature flows that have been extensively studied on surfaces
include Chow-Luo's combinatorial Ricci flow for Thurston's circle pakcings \cite{CL},
Luo's combinatorial Yamabe flow for vertex scaling  \cite{L1},  the combinatorial Calabi flow
for discrete conformal structures \cite{Ge-thesis,Ge,GH1, GX3,X9,ZX} and others,
which were invented and studied independently in the history.
Recently, Chow-Luo's combinatorial Ricci flow for Thurston's circle packings and
Luo's combinatorial Yamabe flow for vertex scalings on surfaces have been unified as the combinatorial Ricci flow
in the framework of discrete conformal structures on polyhedral surfaces \cite{ZGZLYG}.
In this paper, we introduce a fractional combinatorial Calabi flow for discrete conformal structures on polyhedral surfaces,
which unifies and generalizes the combinatorial Ricci flow and combinatorial Calabi flow for discrete conformal structures on surfaces.

Suppose $(M, \mathcal{T})$ is a triangulated connected closed surface with the triangulation $\mathcal{T}=\{V, E, F\}$, where
$V, E, F$ represent the sets of vertices, edges and faces respectively.
For simplicity, we set $|V|=N$ and
use $i$, $\{ij\}$, $\triangle ijk$ to denote the elements in $V, E, F$ respectively.
$\varepsilon: V\rightarrow \{0,1\}$ and
$\eta: E\rightarrow \mathbb{R}$ are two weights defined on the sets of vertices and edges respectively.
$(M, V)$  ($(M, V, \varepsilon)$ respectively) is called as a marked surface (weighted marked surface respectively).
 The following unified notion of discrete conformality was proposed by Glickenstein et al.

\begin{definition}[\cite{G3, GT, T, ZGZLYG}]\label{definition of DCS}
A discrete conformal structure on a weighted triangulated surface $(M, \mathcal{T}, \varepsilon, \eta)$ is a
map $f: V\rightarrow \mathbb{R}$ determining a discrete polyhedral metric $l: E\rightarrow (0,+\infty)$ with
\begin{equation}\label{length of Euclid DCS}
\begin{aligned}
l_{ij}=\sqrt{\varepsilon_ie^{2f_i}+\varepsilon_je^{2f_j}+2\eta_{ij}e^{f_i+f_j}}
\end{aligned}
\end{equation}
in the Euclidean background geometry and
\begin{equation}\label{length of hyper DCS}
\begin{aligned}
 l_{ij}=\cosh^{-1}\left(\sqrt{(1+\varepsilon_ie^{2f_i})(1+\varepsilon_je^{2f_j})}+\eta_{ij}e^{f_i+f_j}\right)
 \end{aligned}
\end{equation}
in the hyperbolic background geometry.
\end{definition}

To determine a discrete polyhedral metric on $(M, \mathcal{T})$, the map $l:E\rightarrow (0,+\infty)$ should satisfy
the triangle inequalities for every face $\triangle ijk\in F$.
The discrete conformal structure in Definition \ref{definition of DCS} unifies and generalizes the
existing special types of discrete conformal structures on surfaces, including
the tangential circle packings ($\varepsilon\equiv1, \eta\equiv1$), Thurston's circle packings ($\varepsilon\equiv1, \eta\in [0,1]$),
inversive distance circle packings ($\varepsilon\equiv1, \eta\in (-1,+\infty)$),
the vertex scaling ($\varepsilon\equiv0, \eta\in (0,+\infty)$) and others.
The discrete conformal structures in Definition \ref{definition of DCS}
could be defined for more general settings, including $\varepsilon_i=-1$ for
some vertices $i\in V$ and the spherical background geometry.
Please refer to \cite{GT, ZGZLYG, X9} for more information on this.
In this paper, we focus on the case that $\varepsilon: V\rightarrow \{0,1\}$ and the Euclidean and hyperbolic background geometry.

Set
\begin{equation}\label{u Euclidean}
\begin{aligned}
u_i=f_i
\end{aligned}
\end{equation}
for any vertex $i\in V$ in the Euclidean background geometry and
\begin{equation}\label{u hyperbolic}
\begin{aligned}
u_i=\left\{
      \begin{array}{ll}
        f_i, & \hbox{ if $\varepsilon_i=0$,} \\
        \frac{1}{2}\log \left|\frac{\sqrt{1+e^{2f_i}}-1}{\sqrt{1+e^{2f_i}}+1}\right|, & \hbox{ if $\varepsilon_i= 1$,}
      \end{array}
    \right.
\end{aligned}
\end{equation}
in the hyperbolic background geometry.
For simplicity, we also call $u$ defined by (\ref{u Euclidean}) and (\ref{u hyperbolic})
as a Euclidean and hyperbolic discrete conformal structure respectively.
The combinatorial Ricci flow
for discrete conformal structure on surfaces \cite{CL,L1,ZGZLYG} is defined to be
\begin{equation}\label{CRF}
\begin{aligned}
\frac{du_i}{dt}=-(K-\overline{K})_i,
\end{aligned}
\end{equation}
where $K: V\rightarrow (-\infty, 2\pi)$ is the combinatorial curvature with $K_i$ defined as $2\pi$ less the cone angle at $i\in V$,
 $\overline{K}: V\rightarrow (-\infty, 2\pi)$ is a fixed function representing a target combinatorial curvature with $\sum_{i\in V}\overline{K}_i=2\pi\chi(M)$
in the Euclidean background geometry and $\sum_{i\in V}\overline{K}_i>2\pi\chi(M)$
in the hyperbolic background geometry.
The combinatorial Ricci flow (\ref{CRF}) introduced in \cite{ZGZLYG} unifies and generalizes
Chow-Luo's combinatorial Ricci flow for Thurston's
circle packings \cite{CL} and Luo's combinatorial Yamabe flow for vertex scaling \cite{L1} on surfaces.
The combinatorial Calabi flow for discrete conformal structures on surfaces \cite{Ge-thesis,Ge,GH1,GX3,X9,ZX} is defined to be
\begin{equation}\label{CCF}
\begin{aligned}
\frac{du_i}{dt}=\Delta(K-\overline{K})_i,
\end{aligned}
\end{equation}
where $\Delta=-L=-(\frac{\partial K}{\partial u})$ is the discrete Laplace operator
for the discrete conformal structures in Definition \ref{definition of DCS}.
The discrete Laplace operator $\Delta$ is proved \cite{X9} to be negative definite (negative semi-definite with rank $N-1$ in the Euclidean background geometry) under the structure condition
\begin{equation}\label{structure condition}
\begin{aligned}
\varepsilon_s \varepsilon_t +\eta_{st}>&0, \ \ \forall \{st\}\in E,\\
\varepsilon_q\eta_{st}+\eta_{qs}\eta_{qt}\geq &0, \ \ \forall \{q s t\}\in F.
\end{aligned}
\end{equation}
See also \cite{BPS, L1, XZ1} for the special case of vertex scaling and \cite{Guo, L4,X1,X3,Z1} for the special case of inversive distance circle packings.
In the case of vertex scaling,
Gu-Luo-Wu \cite{GLW}, Luo-Sun-Wu \cite{LSW} and Wu-Zhu \cite{WZ} recently introduced the following type combinatorial curvature flow
\begin{equation}\label{inversive CCF}
\begin{aligned}
\frac{du_i}{dt}=\Delta^{-1}(K_0-\overline{K})_i
\end{aligned}
\end{equation}
to study the convergence of discrete uniformizaiton conformal factor to the smooth uniformization conformal factor, where
$K_0$ is the initial combinatorial curvature.
Note that the definitions of the three different combinatorial curvature flows (\ref{CRF}), (\ref{CCF}) and  (\ref{inversive CCF})
involve three linear operators $-Id$, $\Delta$ and $\Delta^{-1}$ acting on $K-\overline{K}$ or $K_0-\overline{K}$,
which can be written in a unified form $\Delta^n=-L^n=-(\frac{\partial K}{\partial u})^n$ with $n=0, 1, -1$ respectively.
This motivates us to consider the following fractional combinatorial Laplace operator $\Delta^s$ for $s\in \mathbb{R}$.

As the matrix $L=(\frac{\partial K}{\partial u})$ is symmetric and
positive definite \cite{X9} under the structure condition (\ref{structure condition}),
by the Gram-Schmidt orthonormalization,
there exists an orthonormal matrix $P\in O(N)$ such that
\begin{equation*}
\begin{aligned}
L=\left(\frac{\partial K}{\partial u}\right)=P^T\cdot \text{diag}\{\lambda_1, \cdots, \lambda_n\}\cdot P,
\end{aligned}
\end{equation*}
where $\lambda_1\geq\cdots\geq\lambda_n\geq 0$ are nonnegative eigenvalues of  $L=(\frac{\partial K}{\partial u})$.
For $s\in \mathbb{R}$,
the fractional discrete Laplace operator $\Delta^s$ of order $2s$ \cite{Bell}  is defined to be the matrix
\begin{equation}\label{fractional discrete Laplace}
\begin{aligned}
\Delta^s=-L^s=-P^T\cdot \text{diag}\{\lambda_1^s, \cdots, \lambda_n^s\}\cdot P,
\end{aligned}
\end{equation}
where $0^s$ is set to be $0$ for any $s\in \mathbb{R}$.

The fractional Laplace operators have been one of the most studied research topics in the present century since the work of
Caffarelli-Silvestre \cite{CS} and Silvestre \cite{Si}. It has lots of applications in harmonic analysis, fractional calculus, functional analysis and probability. Especially, it gives a good description of the approximation of discrete jump models to
continuous jump models in random walk \cite{MK}.
Using the geometric fractional Laplace operator, fractional curvature flows have been studied in Riemannian geometry.
A typical example is Jin-Xiong's fractional Yamabe flow \cite{JX} used to study the fractional Yamabe problem.
The fractional discrete Laplace operator in (\ref{fractional discrete Laplace}) is
a discrete analogue of the classical fractional Laplace operator, which has been  extensively
studied in complex networks. Please refer to \cite{RM} and the references therein.
Motivated by the three combinatorial curvature flows (\ref{CRF}), (\ref{CCF}) and  (\ref{inversive CCF}),
we introduce the following fractional combinatorial Calabi flow for discrete conformal structures on surfaces.

\begin{definition}\label{FCCF definition}
Suppose $(M, \mathcal{T}, \varepsilon, \eta)$ is a weighted triangulated connected closed surface with the weights
$\varepsilon: V\rightarrow \{0,1\}$ and $\eta: E\rightarrow \mathbb{R}$ satisfying the structure condition (\ref{structure condition}).
$s\in \mathbb{R}$ is a constant.
The fractional combinatorial Calabi flow of order $s$ for the discrete conformal structures
on $(M, \mathcal{T}, \varepsilon, \eta)$ is defined to be
\begin{equation}\label{FCCF equation}
\begin{aligned}
\frac{du_i}{dt}=\Delta^{s}(K-\overline{K})_i,
\end{aligned}
\end{equation}
where $\Delta^{s}$ is the fractional discrete Laplace operator defined by (\ref{fractional discrete Laplace}).
\end{definition}

\begin{remark}
The fractional combinatorial Calabi flow (\ref{FCCF equation}) unifies and generalizes
Chow-Luo's combinatorial Ricci flow for Thurston's circle packings,
Luo's combinatorial Yamabe flow for vertex scaling
and the combinatorial Calabi flow for discrete conformal structures on surfaces.
Specially, in the cases of $s=0$ and $1$, the fractional combinatorial Calabi flow (\ref{FCCF equation})
is reduced to the combinatorial Ricci flow (\ref{CRF}) and combinatorial Calabi flow (\ref{CCF}) respectively,
which are gradient flows.
Note that for generic $s\in \mathbb{R}$, the fractional combinatorial
Calabi flow (\ref{FCCF equation}) is not a gradient flow.
In the case of $s=-1$, the fractional combinatorial Calabi flow (\ref{FCCF equation}) is
slightly different from the combinatorial curvature flow (\ref{inversive CCF}), where $K$ depends on $t$ in (\ref{FCCF equation}) and $K_0$ is fixed in (\ref{inversive CCF}).
Note that the fractional discrete Laplace operator $\Delta^{s}$
is a nonlocal operator in general (exceptional cases include $s\in \mathbb{Z}_{\geq0}$),
because the eigenvalues $\lambda_1, \cdots, \lambda_n$ globally depend
on the elements of the matrix $L=(\frac{\partial K}{\partial u})$.
This implies that the fractional combinatorial Calabi flow (\ref{FCCF equation}) is a nonlocal combinatorial curvature flow in general.
The fractional combinatorial Calabi flow (\ref{FCCF equation})
is different from the combinatorial $p$-th Calabi flow defined by discrete $p$-Laplace operator in \cite{LZ, FLZ}.
Motivated by Definition \ref{FCCF definition},
we further introduce a fractional combinatorial Calabi flow for decorated and hyper-ideal hyperbolic polyhedral metrics
on $3$-dimensional manifolds in \cite{WX}, where the basic properties of the flow are also established.
Using a fractional discrete Laplace operator for tangential sphere packing metrics on $3$-dimensional manifolds,
a similar fractional combinatorial curvature flow for $s\geq 0$ was previously introduced in \cite{GX1}.
\end{remark}

As the eigenvalues $\lambda_1, \cdots, \lambda_n$ of  $L=(\frac{\partial K}{\partial u})$
are Lipschitz functions of the discrete conformal structures,
the short time existence for the solution of  fractional combinatorial Calabi flow (\ref{FCCF equation})
follows by the standard theory in ordinary differential equations.
For the fractional combinatorial Calabi flow (\ref{FCCF equation}),
we further have the following result.

\begin{theorem}\label{main theorem local convergence of FCCF}
Suppose $(M, \mathcal{T}, \varepsilon, \eta)$ is a weighted triangulated connected closed surface with the weights
$\varepsilon: V\rightarrow \{0,1\}$ and $\eta: E\rightarrow \mathbb{R}$ satisfying the structure condition (\ref{structure condition}).
If there exists a discrete conformal structure $\overline{u}$ with combinatorial curvature $\overline{K}$,
then for any $s\in \mathbb{R}$,
there exists a positive constant $\delta>0$ such that if $||u(0)-\overline{u}||<\delta$
($\sum_{i\in V}u_i(0)=\sum_{i\in V}\overline{u}$ additionally in the case of Euclidean background geometry),
 the solution $u(t)$ of the fractional combinatorial Calabi flow (\ref{FCCF equation})
 exists for all time and converges exponentially fast to $\overline{u}$.
\end{theorem}

For generic initial discrete conformal structures in Definition \ref{definition of DCS},
the solution of fractional combinatorial Calabi flow (\ref{FCCF equation}) on $(M, \mathcal{T}, \varepsilon, \eta)$
may develop singularities, which correspond to the triangles degenerate along the flow (\ref{FCCF equation})
or the conformal factors tend to infinity.
In the case of $s=0$, the second author \cite{X9} proved the longtime existence and
global convergence for the solution of fractional combinatorial Calabi flow (\ref{FCCF equation})
by extending the flow through singularities.
For $s\neq 0$, we do not have any unified approach to handle the singularities
along the fractional combinatorial Calabi flow (\ref{FCCF equation})
for the discrete conformal structures in Definition \ref{definition of DCS}.
However, for Thurston's circle packings on surfaces,
we show that the singularities never develop along (\ref{FCCF equation})
and prove the following longtime existence and global convergence for
the solution of fractional combinatorial Calabi flow  (\ref{FCCF equation}).
\begin{theorem}\label{main theorem conv of FCCF Thurston's CP}
Suppose $(M, \mathcal{T}, \varepsilon, \eta)$ is a weighted triangulated connected closed surface with the weights
$\varepsilon\equiv 1$ and $\eta: E\rightarrow (-1, 1]$ satisfying the structure condition (\ref{structure condition}).
\begin{description}
  \item[(a)] In the Euclidean background geometry, if there exists a discrete conformal structure $\overline{u}\in \mathbb{R}^N$ with combinatorial curvature $\overline{K}$,  then for any $s\in \mathbb{R}$ and any initial value $u(0)\in \mathbb{R}^N$ with
      $\sum_{i\in V}u_i(0)=\sum_{i\in V}\overline{u}_i$,
      the solution of fractional combinatorial Calabi flow  (\ref{FCCF equation})
     exists for all time and converges exponentially fast to $\overline{u}$.
  \item[(b)] In the hyperbolic background geometry, if there exists a discrete conformal structure
  $\overline{u}\in \mathbb{R}^N_{<0}$ with combinatorial curvature $\overline{K}$,
  then for any $s\in \mathbb{R}$ and initial value $u(0)\in \mathbb{R}^N_{<0}$,
      the solution of fractional combinatorial Calabi flow  (\ref{FCCF equation}) exists for all time and converges exponentially fast to $\overline{u}$.
\end{description}
\end{theorem}

\begin{remark}
If $s=0$ and $\eta\in [0,1]$,
the result in Theorem \ref{main theorem conv of FCCF Thurston's CP} is obtained by Chow-Luo \cite{CL}
for combinatorial Ricci flow of Thurston's circle packings.
If $s=1$ and $\eta\in [0,1]$, the result in Theorem \ref{main theorem conv of FCCF Thurston's CP} is reduced to the results
obtained in \cite{Ge-thesis, Ge,GX3,GH1}  for combinatorial Calabi flow of Thurston's circle packings.
Following Thurston's arguments in \cite{T1} word by word,
one can replace the condition on the existence of $\overline{u}$ with combinatorial curvature $\overline{K}$
in Theorem \ref{main theorem conv of FCCF Thurston's CP}
by some linear equalities and inequalities on $\overline{K}$, which characterize the image of the curvature map $K$
for Thurston's circle packings with $\eta: E\rightarrow (-1, 1]$ satisfying the structure condition (\ref{structure condition}).
One can also refer to \cite{CL, GHZ,X1arxiv-v1} for this.
\end{remark}

For vertex scaling,
the fractional combinatorial Calabi flow  (\ref{FCCF equation}) on $(M, \mathcal{T}, \varepsilon, \eta)$ may develop singularities.
In this case,
we do surgery on the fractional combinatorial Calabi flow (\ref{FCCF equation}) by edge flipping under the Delaunay condition
to handle the potential singularities along  (\ref{FCCF equation}),
which was first introduced in \cite{GGLSW,GLSW} for combinatorial Yamabe flow.
Here we give a brief description of the surgery in the Euclidean background geometry.
For a piecewise linear metric (PL metric for short in the following) defined on $(M, \mathcal{T})$, it is said to satisfy the Delaunay condition
if for every edge $\{ij\}\in E$ we have $\theta_{k}^{ij}+\theta_{l}^{ij}\leq \pi$,
where $\theta_{k}^{ij}, \theta_{l}^{ij}$ are inner angles
facing the edge $\{ij\}\in E$ in the triangles $\triangle ijk$ and $\triangle ijl$ respectively.
Along the fractional combinatorial Calabi flow  (\ref{FCCF equation}) on $(M, \mathcal{T}, \varepsilon, \eta)$,
suppose the Delaunay condition is satisfied for $t\in [0, T]$ and there exists an edge $\{ij\}\in E$ and a constant $\epsilon>0$ such that
for any $t\in (T, T+\epsilon)$, we have $\theta_{k}^{ij}+\theta_{l}^{ij}> \pi$.
In this case, we replace the edge $\{ij\}\in E$ by a new edge $\{kl\}$ to get a new triangulation $\mathcal{T}'$ at the time $t=T$
and then evolve the fractional combinatorial Calabi flow  (\ref{FCCF equation})
with the PL metric on $(M, \mathcal{T})$ at $t=T$ as the initial metric on
$(M, \mathcal{T}')$. This process is called \textbf{\emph{surgery by edge flipping}} under the Delaunay condition.
The surgery by edge flipping could also be defined for piecewise hyperbolic metrics (PH metrics for short in the following) on surfaces
under the hyperbolic Delaunay condition, which is defined to be
$\theta_{k}^{ij}+\theta_{l}^{ij}\leq \theta_{i}^{jk}+\theta_{i}^{jl}+\theta_{j}^{ik}+\theta_{j}^{il}$
for adjacent triangles $\triangle ijk, \triangle ijl\in F$.
With the help of discrete conformal theory established in \cite{GLSW, GGLSW},
we prove the following result for the
fractional combinatorial Calabi flow with surgery in the case of vertex scaling.

\begin{theorem}\label{main theorem conv of FCCF vertex scaling}
Suppose $(M, V)$ is a connected closed marked surface.
\begin{description}
  \item[(a)] In the Euclidean background geometry, if $\overline{K}: V\rightarrow (-\infty, 2\pi)$ satisfies $\sum_{i\in V} \overline{K}_i=2\pi\chi(M)$, then for any $s\in \mathbb{R}$ and any initial PL metrics on $(M, V)$,
      the solution of fractional combinatorial Calabi flow with surgery on $(M, V)$
      for vertex scaling exists for all time and converges exponentially fast.
  \item[(b)] In the hyperbolic background geometry, if $\overline{K}: V\rightarrow (-\infty, 2\pi)$ satisfies $\sum_{i\in V} \overline{K}_i>2\pi\chi(M)$, then for any $s\in \mathbb{R}$ and any initial PH metrics on $(M, V)$,
      the solution of fractional combinatorial Calabi flow with surgery on $(M, V)$
      for vertex scaling exists for all time and converges exponentially fast.
\end{description}
\end{theorem}

\begin{remark}
If $s=0$, the result in Theorem \ref{main theorem conv of FCCF vertex scaling} was proved by Gu-Luo-Sun-Wu \cite{GLSW}
in the Euclidean background geometry
and by Gu-Guo-Luo-Sun-Wu \cite{GGLSW} in the hyperbolic background geometry respectively
for combinatorial Yamabe flow with surgery.
If $s=1$, the result in Theorem \ref{main theorem conv of FCCF vertex scaling}
was proved by Zhu-Xu \cite{ZX} for combinatorial Calabi flow with surgery.
As the fractional combinatorial Calabi flow (\ref{FCCF equation}) is not a gradient flow in general,
the result proved by the first author in \cite{Wu} can not be applied directly to prove that the number of surgeries along the
fractional combinatorial Calabi flow (\ref{FCCF equation}) is finite.
It is conceived that for any $s\in \mathbb{R}$, this is true for fractional combinatorial Calabi flow with surgery.
\end{remark}

For discrete conformal structures in Definition \ref{definition of DCS},
one can introduce the notion of weighed Delaunay condition.
For a PL or PH metric on $(M, \mathcal{T}, \varepsilon,\eta)$
generated by discrete conformal structures in Definition \ref{definition of DCS},
it is said to satisfy the \textbf{\textit{weighted Delaunay condition}}
if $\frac{\partial K_i}{\partial u_j}\leq 0$ for every edge $\{ij\}\in E$.
In the case of $\varepsilon\equiv 0$, which corresponds to vertex scaling,
the weighted Delaunay condition is equivalent to the standard Delaunay condition.
Please refer to \cite{BPS,BS, G2a, G3,G4,GT,Lei,WZ, X9, XZ1,XZ2} for more information on Delaunay condition and
weighted Delaunay condition.
For the fractional combinatorial Calabi flow  (\ref{FCCF equation}) of discrete conformal structures in Definition \ref{definition of DCS},
one can also do surgery along the flow (\ref{FCCF equation}) by edge flipping under the weighted Delaunay condition.
We have the following conjecture on the global convergence of the solution of
fractional combinatorial Calabi flow with surgery for discrete conformal structures in Definition \ref{definition of DCS}.

\begin{conjecture}\label{conjecture}
Suppose $(M, V, \varepsilon)$ is a marked weighted connected closed surface with
$\varepsilon: V\rightarrow \{0,1\}$ and there exists a PL or PH metric generated by a discrete conformal structure
in Definition \ref{definition of DCS} with combinatorial curvature $\overline{K}$.
For any $s\in \mathbb{R}$ and any initial PL or PH metric on $(M, V, \varepsilon)$ generated by discrete conformal structures
in Definition \ref{definition of DCS},  the solution of
fractional combinatorial Calabi flow with surgery exists for all time and converges exponentially fast
after finite number of surgeries.
\end{conjecture}

In Conjecture \ref{conjecture},
the PL or PH metric generated by some discrete conformal structure in Definition \ref{definition of DCS}
with combinatorial curvature $\overline{K}$ does not depend on the triangulations of $(M, V, \varepsilon)$.
Conjecture \ref{conjecture} is a generalization of  Conjecture 3 in \cite{X9} proposed by the second author
on the convergence of combinatorial Ricci flow with surgery
and combinatorial Calabi flow with surgery for discrete conformal structures in Definition \ref{definition of DCS}.
Note that for Thurston's circle packings on a triangulated surface $(M, \mathcal{T})$, the weighted Delaunay condition
is automatically satisfied, so we do not need to do surgery on the fractional combinatorial Calabi flow (\ref{FCCF equation}).
For vertex scaling on a triangulated surface $(M, \mathcal{T})$, the weighted Delaunay condition is equivalent to the standard
Delaunay condition. This implies that the surgery by edge flipping under the weighted Delaunay condition for vertex scaling
is the same as the surgery by edge flipping under the standard Delaunay condition.
Therefore, Theorem \ref{main theorem conv of FCCF Thurston's CP}
and Theorem \ref{main theorem conv of FCCF vertex scaling} provides
strong supports for Conjecture \ref{conjecture}.

The paper is organized as follows. In Section \ref{section 2}, we study the basic properties of the fractional
combinatorial Calabi flow for discrete conformal structures on surfaces in Definition \ref{definition of DCS}
and prove Theorem \ref{main theorem local convergence of FCCF}.
In Section \ref{section 3}, we study the fractional
combinatorial Calabi flow for Thurston's circle packings on surfaces
and prove Theorem \ref{main theorem conv of FCCF Thurston's CP}.
In Section \ref{section 4}, we study the fractional
combinatorial Calabi flow for vertex scaling on surfaces and prove Theorem \ref{main theorem conv of FCCF vertex scaling}.
\\
\\
\textbf{Acknowledgements}\\[8pt]
The authors thank Professor Feng Luo, Dr. Yanwen Luo,  Dr.  Wai Yeung Lam,
Xiaoping Zhu, Professor Linlin Sun and Professor Changsong Deng and for communications.
The research of the second author is supported by the Fundamental Research Funds for the Central Universities under
grant no. 2042020kf0199.

\section{Basic properties of fractional combinatorial Calabi flow for discrete conformal structures on surfaces}\label{section 2}
Recall the following result on the matrix $L=(\frac{\partial K}{\partial u})$ for the discrete conformal structures in Definition \ref{definition of DCS}, which generalizes the results obtained in \cite{BPS, Guo, L1, L4,X1,X3,XZ1,Z1}.

\begin{theorem}[\cite{X9}]\label{positve definiteness of L}
Suppose $(M, \mathcal{T}, \varepsilon, \eta)$ is a weighted triangulated surface with
the weights $\varepsilon: V\rightarrow \{0, 1\}$ and $\eta: E\rightarrow \mathbb{R}$
satisfying the structure condition (\ref{structure condition}).
\begin{description}
  \item[(a)] In the Euclidean background geometry,
  the matrix $L=(\frac{\partial K}{\partial u})$
is symmetric and positive semi-definite with rank $N-1$ and kernel $\{t\mathbf{1}=t(1, \cdots, 1)^T\in \mathbb{R}^N|t\in \mathbb{R}\}$
for all nondegenerate Euclidean discrete conformal structures on $(M, \mathcal{T}, \varepsilon, \eta)$.
  \item[(b)] In the hyperbolic background geometry,
  the matrix $L=(\frac{\partial K}{\partial u})$
is symmetric and strictly positive definite for all
nondegenerate hyperbolic discrete conformal structures on $(M, \mathcal{T}, \varepsilon, \eta)$.
\end{description}
\end{theorem}
As a direct consequence, we have the following results.
\begin{lemma}\label{positive defi of Ls}
Suppose $(M, \mathcal{T}, \varepsilon, \eta)$ is a weighted triangulated surface with
the weights $\varepsilon: V\rightarrow \{0, 1\}$ and $\eta: E\rightarrow \mathbb{R}$
satisfying the structure condition (\ref{structure condition}).
\begin{description}
  \item[(a)] In the Euclidean background geometry, for any $s\in \mathbb{R}$, the matrix $L^s$ is symmetric and
  positive semi-definite with rank $N-1$
  and kernel $\{t\mathbf{1}=t(1, \cdots, 1)^T\in \mathbb{R}^N|t\in \mathbb{R}\}$
for all nondegenerate Euclidean discrete conformal structures on $(M, \mathcal{T}, \varepsilon, \eta)$.
  \item[(b)] In the hyperbolic background geometry, for any $s\in \mathbb{R}$,
  the matrix $L^s$
is symmetric and strictly positive definite for all
nondegenerate hyperbolic discrete conformal structures on $(M, \mathcal{T}, \varepsilon, \eta)$.
\end{description}
\end{lemma}
\begin{remark}
Note that for a nondegenerate Euclidean discrete conformal structure on $(M, \mathcal{T}, \varepsilon, \eta)$, we have
$\sum_{i\in V}(K-\overline{K})_i=2\pi\chi(M)-2\pi\chi(M)=0$, which implies
$K-\overline{K}\in Ker(L^s)^{\perp}=Im(L^s)=\mathbf{1}^{\perp}$.
By Lemma \ref{positive defi of Ls}, restricted to $Ker(L^s)^{\perp}=Im(L^s)=\mathbf{1}^{\perp}$,
$\Delta^s$ is a nonsingular linear operator.
Specially, $\Delta^0|_{\mathbf{1}^{\perp}}=-Id|_{\mathbf{1}^{\perp}}$ and $\Delta^1|_{\mathbf{1}^{\perp}}=-L|_{\mathbf{1}^{\perp}}$.
Combining this with $K-\overline{K}\in \mathbf{1}^{\perp}$
implies that the $0$-order Euclidean fractional combinatorial Calabi flow (\ref{FCCF equation})
is the Euclidean combinatorial Ricci flow (\ref{CRF})
and the $1$-order Euclidean fractional combinatorial Calabi flow (\ref{FCCF equation})
is the Euclidean combinatorial Calabi flow (\ref{CCF}).
Parallelling results hold for the hyperbolic background geometry.
Therefore, the fractional combinatorial Calabi flow (\ref{FCCF equation}) unifies and generalizes
the  combinatorial Ricci flow (\ref{CRF}) and the combinatorial Calabi flow (\ref{CCF}).
\end{remark}

\begin{remark}
In the Euclidean background geometry, if $\overline{K}_i=K_{av}=\frac{2\pi\chi(M)}{N}$ for every $i\in V$,
the fractional combinatorial Calabi flow (\ref{FCCF equation}) is equivalent to
$$\frac{du_i}{dt}=\Delta^s K_i$$
by Lemma \ref{positive defi of Ls} (a), which is the standard form of fractional combinatorial Calabi flow.
As the results in this paper could be proved for prescribed combinatorial curvature $\overline{K}$, we introduce
the fractional combinatorial Calabi flow in the form (\ref{FCCF equation}) to simplify the statements of the results.
\end{remark}

By Lemma \ref{positive defi of Ls} (a), we further have the
following property for the solution of the Euclidean fractional combinatorial
Calabi flow (\ref{FCCF equation}).

\begin{lemma}\label{lemma invariance of sum ui}
Suppose $(M, \mathcal{T}, \varepsilon, \eta)$ is a weighted triangulated surface with
the weights $\varepsilon: V\rightarrow \{0, 1\}$ and $\eta: E\rightarrow \mathbb{R}$
satisfying the structure condition (\ref{structure condition}).
$u(t)$ is a solution of the Euclidean fractional combinatorial Calabi flow (\ref{FCCF equation}).
Then $\sum_{i=1}^{N}u_i(t)$ is invariant along the Euclidean fractional combinatorial Calabi flow (\ref{FCCF equation}).
\end{lemma}
\proof
By direct calculations, we have
\begin{equation*}
\begin{aligned}
\frac{d}{dt}\left(\sum_{i=1}^{N}u_i(t)\right)=-\sum_{i=1}^{N}L^s(K-\overline{K})_i=-\mathbf{1}^TL^s(K-\overline{K})=0
\end{aligned}
\end{equation*}
by Lemma \ref{positive defi of Ls} (a).
\qed

Lemma \ref{lemma invariance of sum ui} implies that  the solution of
Euclidean fractional combinatorial Calabi flow (\ref{FCCF equation})
stays in the hyperplane $\{u\in \mathbb{R}^N|\sum_{i=1}^Nu_{i}=\sum_{i=1}^Nu_{i}(0)\}$.
Without loss of generality, we will assume $u(0)$ is in the hyperplane $\Sigma_0:=\{u\in \mathbb{R}^N|\sum_{i=1}^Nu_{i}=0\}$
for the Euclidean background geometry in the following.\\

\textbf{Proof of Theorem \ref{main theorem local convergence of FCCF}:}
Set $\Gamma_i(u)=\Delta^s(K-\overline{K})_i$.
By assumption, $\overline{u}$ is a equilibrium point of the system $\frac{du}{dt}=\Gamma(u)$.
Furthermore, $D\Gamma|_{u=\overline{u}}=-L^{s+1}$.
In the case of hyperbolic background geometry, $D\Gamma$ is negative definite at $\overline{u}$ by Lemma \ref{positive defi of Ls} (b).
Then the result in Theorem \ref{main theorem local convergence of FCCF} for the hyperbolic background geometry
is a direct application of the Lyapunov stability theorem (\cite{P} Chapter 5).
In the case of Euclidean background geometry, $D\Gamma|_{u=\overline{u}}=-L^{s+1}$ is negative semi-definite
with kernel $\{t\mathbf{1}\in \mathbb{R}^N|t\in \mathbb{R}\}$ by Lemma \ref{positive defi of Ls} (a),
which is perpendicular to the hyperplane $\Sigma_0$.
Restricted to $\Sigma_0$, then the following of the proof
for Theorem \ref{main theorem local convergence of FCCF} in the Euclidean background geometry
is also a direct application of the Lyapunov stability theorem.
\qed

\section{Fractional combinatorial Calabi flow for Thurston's circle packings on surfaces}\label{section 3}
Thurston's circle packing \cite{T1} is a special type of discrete conformal structures on polyhedral surfaces
with $\varepsilon\equiv 1$ and $\eta\in [0, 1]$ in Definition \ref{definition of DCS}.
Motivated by \cite{X1,Z1,Z2}, we consider a generalization of Thurston's original definition,
which corresponds to $\varepsilon\equiv 1$ and $\eta\in (-1, 1]$ in Definition \ref{definition of DCS}.

\subsection{The case of Thurston's Euclidean circle packings}
In the Euclidean background geometry, the edge length defined by Thurston's circle packing is given by
\begin{equation}\label{Euclidean Thurs CP metric}
  \begin{aligned}
  l_{ij}=\sqrt{e^{2u_i}+e^{2u_j}+2\eta_{ij}e^{u_i+u_j}}.
  \end{aligned}
\end{equation}
We have the following result on the triangle inequality for the edge length defined by (\ref{Euclidean Thurs CP metric}).
\begin{lemma}[\cite{X1} Corollary 2.2]\label{admiss space for Eucl Thurs CP}
Suppose $\eta: E\rightarrow (-1, 1]$ satisfies the structure condition (\ref{structure condition}).
Then for any $u\in \mathbb{R}^N$, the triangle inequalities are satisfied for every triangle
$\triangle ijk\in F$.
\end{lemma}

Combining Theorem \ref{positve definiteness of L} and Lemma \ref{admiss space for Eucl Thurs CP},
the function
$F(u)=\int^u_{0}\sum_{i=1}^N(K_i-\overline{K}_i)du_i$
is a well-defined smooth convex function defined on $\mathbb{R}^N$ with $\nabla_u F=K-\overline{K}$ and $Hess_u F=L$.
If there exists $\overline{u}\in \Sigma_0$ with $K(\overline{u})=\overline{K}$, then
$\nabla_u F(\overline{u})=0$, which implies $F(\overline{u})=\min_{u\in \Sigma_0}F(u)$
and $\lim_{u\rightarrow \infty, u\in \Sigma_0}F(u)=+\infty$ by the convexity of $F(u)$.\\

\textbf{Proof of Theorem \ref{main theorem conv of FCCF Thurston's CP} (a):}
Along the fractional combinatorial Calabi flow (\ref{FCCF equation}), we have
\begin{equation*}
\begin{aligned}
\frac{d}{dt}F(u(t))=\sum_{i=1}^N\nabla_{u_i}F\cdot \frac{du_i}{dt}=-(K-\overline{K})^T\cdot L^s\cdot (K-\overline{K})\leq 0
\end{aligned}
\end{equation*}
by Lemma \ref{positive defi of Ls} (a), which implies $F(u(t))\leq F(u(0))$ along the fractional combinatorial Calabi flow (\ref{FCCF equation}).
Combining Lemma \ref{lemma invariance of sum ui} and $\lim_{u\rightarrow \infty, u\in \Sigma_0}F(u)=+\infty$,
the solution $u(t)$ of Euclidean fractional combinatorial Calabi flow (\ref{FCCF equation}) stays in a compact subset $\Omega$ of $\Sigma_0$,
which further implies the solution $u(t)$ of  (\ref{FCCF equation}) exists for all time.

As the solution $u(t)$ of the Euclidean fractional combinatorial Calabi flow (\ref{FCCF equation})
stays in a compact subset $\Omega\subset\subset\Sigma_0$
and $L^{s+1}$ is strictly positive definite on $\Sigma_0$ by Lemma \ref{positive defi of Ls},
the first nonzero eigenvalue of $L^{s+1}$, which is a continuous function of
the discrete conformal structures, has a positive lower bound $\lambda$
along the Euclidean fractional combinatorial Calabi flow (\ref{FCCF equation}).
Therefore, for the combinatorial Calabi energy $\overline{\mathcal{C}}(t):=\sum_{i=1}^N(K_i-\overline{K}_i)^2$, we have
\begin{equation*}
\begin{aligned}
\frac{d}{dt}\overline{\mathcal{C}}(t)
=2\sum_{i=1}^N(K_i-\overline{K}_i)\frac{dK_i}{dt}
=-2(K-\overline{K})^S\cdot L^{s+1}\cdot (K-\overline{K})
\leq -2\lambda \overline{\mathcal{C}}(t),
\end{aligned}
\end{equation*}
which implies $\overline{\mathcal{C}}(t)\leq e^{-2\lambda t}\overline{\mathcal{C}}(0)$.
As $K|_{\Sigma_0}$ is a diffeomorphism from $\Sigma_0$ to $K(\Sigma_0)$ by Theorem \ref{positve definiteness of L},
this further implies that the solution $u(t)$ of Euclidean fractional combinatorial Calabi flow (\ref{FCCF equation})
converges exponentially fast to $\overline{u}$.
\qed

\subsection{The case of Thurston's hyperbolic circle packings}
In the hyperbolic background geometry, set $e^{f_i}=\sinh r_i$
in (\ref{length of hyper DCS}).
Then the edge length $l_{ij}$ defined by Thurston's hyperbolic circle packing is given by
 \begin{equation}\label{hyper length in r for Th CP}
\begin{aligned}
 \cosh l_{ij}=\cosh r_i\cosh r_j+\eta_{ij}\sinh r_i\sinh r_j.
 \end{aligned}
\end{equation}
The map $r: V\rightarrow (0, +\infty)$ is called Thurston's hyperbolic circle packing metric.
For simplicity of notations, set
\begin{equation}\label{Ci Si}
\begin{aligned}
C_i=\cosh r_i, S_i=\sinh r_i
 \end{aligned}
\end{equation}
in the following of this subsection.

\begin{lemma}[\cite{X9} Lemma 4.6]\label{equi of lij Thurs hyper CP}
Suppose $\eta_{ij}\in (-1,1]$. For the edge length $l_{ij}$ defined by (\ref{hyper length in r for Th CP}),
there exist $C=C(\eta_{ij})>0$ and $C'=C'(\eta_{ij})>0$ such that
\begin{equation*}
\begin{aligned}
C(C_iC_j+S_iS_j)\leq \cosh l_{ij}\leq C'(C_iC_j+S_iS_j).
 \end{aligned}
\end{equation*}
\end{lemma}
Parallelling to Lemma \ref{admiss space for Eucl Thurs CP} for Thurston's Euclidean circle packing metrics,
we have the following result on triangle inequalities for Thurston's hyperbolic circle packing metrics.

\begin{lemma}[\cite{Z1} Lemma 2.4, \cite{X1} Corollary 3.2]\label{admi space for Thurs Eucl CP}
Suppose $\eta: E\rightarrow (-1, 1]$ satisfies the structure condition (\ref{structure condition}).
For any $r\in \mathbb{R}^N_{>0}$,
the triangle inequalities are satisfied for any triangle
$\triangle ijk\in F$ with edge length defined by (\ref{hyper length in r for Th CP}).
\end{lemma}

By the definition of $u_i$ in (\ref{u hyperbolic}), we have
 \begin{equation}\label{ui ri relationship}
\begin{aligned}
u_i
=\frac{1}{2}\log\frac{\cosh r_i-1}{\cosh r_i+1}
=\log\tanh\frac{r_i}{2}
 \end{aligned}
\end{equation}
for Thurston's hyperbolic circle packing metrics,
which implies $u\in \mathbb{R}^N_{<0}$ by Lemma \ref{admi space for Thurs Eucl CP}.
By Theorem \ref{positve definiteness of L} (b) and Lemma \ref{admi space for Thurs Eucl CP}, the function
\begin{equation*}
\begin{aligned}
F(u)=\int^u_{u_0}\sum_{i=1}^N(K_i-\overline{K}_i)du_i
\end{aligned}
\end{equation*}
is a strictly convex function defined on $\mathbb{R}^N_{<0}$ with $\nabla_uF=K-\overline{K}$.
If there exists $\overline{u}\in \mathbb{R}^N_{<0}$ with $K(\overline{u})=\overline{K}$,
then $\nabla_uF(\overline{u})=0$, which further implies that
$\lim_{u\rightarrow \infty, u\in \mathbb{R}^N_{<0}}F(u)=+\infty$ by the strict convexity of $F(u)$.

\begin{lemma}\label{lower bound of u_i(t) for hyperbolic Thurs CP}
Suppose $(M, \mathcal{T}, \varepsilon, \eta)$ is a weighted triangulated connected closed surface with the weights
$\varepsilon\equiv 1$ and $\eta: E\rightarrow (-1, 1]$ satisfying the structure condition (\ref{structure condition}).
Suppose that there exists  $\overline{u}\in \mathbb{R}^N_{<0}$ with $K(\overline{u})=\overline{K}$
in the hyperbolic background geometry and
$u(t)$ is a solution of the hyperbolic fractional combinatorial Calabi flow (\ref{FCCF equation}).
Then $u_i(t)$ is uniformly bounded from below along the hyperbolic fractional combinatorial Calabi flow (\ref{FCCF equation})
for every $i\in V$.
\end{lemma}
\proof
Along the hyperbolic fractional combinatorial Calabi flow (\ref{FCCF equation}), we have
\begin{equation*}
\begin{aligned}
\frac{d}{dt}F(u(t))=\sum_{i=1}^N\nabla_{u_i}F\cdot \frac{du_i}{dt}=-(K-\overline{K})^T\cdot L^s\cdot (K-\overline{K})\leq 0
\end{aligned}
\end{equation*}
by Lemma \ref{positive defi of Ls} (b), which implies $F(u(t))\leq F(u(0))$
along the hyperbolic fractional combinatorial Calabi flow (\ref{FCCF equation}).
Combining with the fact that $\lim_{u\rightarrow \infty, u\in \mathbb{R}^N_{<0}}F(u)=+\infty$
under the existence of $\overline{u}\in \mathbb{R}^N_{<0}$ with $K(\overline{u})=\overline{K}$,
we have $u(t)$ is uniformly bounded, which implies that
 $u_i(t)$ is uniformly bounded from below along the hyperbolic fractional combinatorial Calabi flow (\ref{FCCF equation})
for every $i\in V$.
\qed

\begin{remark}\label{lower bound of r_i geq R}
By (\ref{ui ri relationship}),
$r_i\in (0, +\infty)$ is a strictly nondecreasing function of $u_i\in (-\infty, 0)$ with
$\lim_{u_i\rightarrow -\infty}r_i=0$ and $\lim_{u_i\rightarrow 0^-}r_i=+\infty$.
Under the existence of Thurston's hyperbolic circle packing with combinatorial curvature $\overline{K}$,
Lemma \ref{lower bound of u_i(t) for hyperbolic Thurs CP} implies that
there exists a constant $R>0$ such that $r_i(t)>R$ for all $i\in V$ along the fractional combinatorial Calabi flow (\ref{FCCF equation}).
We will always assume $r_i\geq R>0$ in the following of this subsection.
By the proof of Lemma \ref{lower bound of u_i(t) for hyperbolic Thurs CP},
one can also get a positive upper bound for $u_i(t)$.
However, $u_i<0$ by  (\ref{ui ri relationship}). Therefore, the positive upper bound of $u_i(t)$ is useless.
\end{remark}

\begin{lemma}\label{lower bound of ch lij-1}
Suppose $\eta_{ij}\in (-1, 1]$ for every edge $\{ij\}\in E$ and there exists a constant $R>0$
such that $r_i>R$ for every $i\in V$.  Then $\cosh l_{ij}-1$ has
a uniform lower bound $C=C(\eta, R)>0$ for every edge $\{ij\}\in E$.
\end{lemma}
\proof
By Cauchy inequality, we have
\begin{equation*}
\begin{aligned}
\cosh l_{ij}
=&\sqrt{(1+\sinh^2 r_i)(1+\sinh^2 r_j)}+\eta_{ij}\sinh r_i \sinh r_j\\
\geq & 1+\sinh r_i\sinh r_j+\eta_{ij}\sinh r_i \sinh r_j\\
=&1+(1+\eta_{ij})\sinh r_i\sinh r_j,
 \end{aligned}
\end{equation*}
which implies
\begin{equation*}
\begin{aligned}
\cosh l_{ij}-1
\geq (1+\eta_{ij})\sinh r_i\sinh r_j\geq \min_{\{ij\}\in E} (1+\eta_{ij}) \sinh^2R>0
 \end{aligned}
\end{equation*}
by $\eta_{ij}\in (-1, 1]$ for every edge $\{ij\}\in E$ and $r_i\geq R$ for every $i\in V$.
We can take $C=\min_{\{ij\}\in E} (1+\eta_{ij}) \sinh^2R$.
\qed

Recall that for Thurston's hyperbolic circle packing metrics,
the matrix $L=(\frac{\partial K}{\partial u})$ could be decomposed as (\cite{GX3} Theorem 3.1)
\begin{equation*}
\begin{aligned}
L=A+B,
 \end{aligned}
\end{equation*}
where $A$ is a diagonal matrix with
$$A_{ii}=\frac{\partial}{\partial u_i}(\sum_{\triangle ijk\in F}Area(\triangle ijk))$$
and $B$ is a symmetric matrix with
\begin{equation}\label{defin of B}
\begin{aligned}
B_{ij}=\left\{
         \begin{array}{ll}
           -(\frac{\partial \theta_{i}^{jk}}{\partial u_j}+\frac{\partial \theta_{i}^{jl}}{\partial u_j}), & \hbox{ if $j\sim i$;} \\
           -\sum_{k\sim i}B_{ik}, & \hbox{ if  $j=i$;} \\
           0, & \hbox{otherwise.}
         \end{array}
       \right.
 \end{aligned}
\end{equation}
Further recall the following formula (\cite{Z1} Lemma 4.1, \cite{X1} Lemma 3.6)
for the derivative of inner angle $\theta_i^{jk}$ in the triangle $\triangle ijk$
\begin{equation}\label{derivative of theta Thurs hyper CP}
\frac{\partial \theta_i^{jk}}{\partial u_j}=\frac{\partial \theta_j^{ik}}{\partial u_i}
=\frac{1}{A_{ijk}\sinh^2 l_{ij}}[C_kS_i^2S_j^2(1-\eta_{ij}^2)+C_iS_iS_j^2S_k\gamma_{jik}+C_jS^2_iS_jS_k\gamma_{ijk}],
\end{equation}
where $A_{ijk}=\sinh l_{ij}\sinh l_{ik}\sin \theta_i^{jk}$ and $\gamma_{ijk}=\eta_{jk}+\eta_{ij}\eta_{ik}$.
By (\ref{derivative of theta Thurs hyper CP}), $\eta\in (-1, 1]$  and the structure condition (\ref{structure condition}),
we have
\begin{equation}\label{theta deriv nonnegative Bij nonpositive}
\frac{\partial \theta_i^{jk}}{\partial u_j}\geq 0, \frac{\partial \theta_i^{jl}}{\partial u_j}\geq 0,
-B_{ij}=\frac{\partial \theta_{i}^{jk}}{\partial u_j}+\frac{\partial \theta_{i}^{jl}}{\partial u_j} \geq 0
\end{equation}
for $j\sim i$.
Note that $\frac{\partial \theta_i^{jk}}{\partial u_j}= 0$ if and only if $\eta_{ij}=1$ and $\eta_{jk}+\eta_{ik}=0$,
which is attainable.
By the following formula obtained by Glickenstein and Thomas (\cite{GT}, Proposition 9)
\begin{equation}\label{GT's formula}
\frac{\partial }{\partial u_i} Area(\triangle ijk)=\frac{\partial \theta_j^{ik}}{\partial u_i}(\cosh l_{ij}-1)+\frac{\partial \theta_k^{ij}}{\partial u_i}(\cosh l_{ik}-1),
\end{equation}
we have
\begin{equation}\label{Aii expression}
\begin{aligned}
A_{ii}=\sum_{j\sim i}(\frac{\partial \theta_i^{jk}}{\partial u_j}+\frac{\partial \theta_i^{jl}}{\partial u_j})(\cosh l_{ij}-1)
=\sum_{j\sim i}(-B_{ij})(\cosh l_{ij}-1)
\geq 0
 \end{aligned}
\end{equation}
by (\ref{theta deriv nonnegative Bij nonpositive}) ,
where $\triangle ijk$ and $\triangle ijl$ are adjacent triangles sharing the common edge $\{ij\}\in E$.
In the case that $r_i>R>0$ for every $i\in V$, we further have the following stronger result on $A_{ii}$.

\begin{lemma}\label{estimate of A-ii}
Suppose $(M, \mathcal{T}, \varepsilon, \eta)$ is a weighted triangulated connected closed surface with the weights
$\varepsilon\equiv 1$ and $\eta: E\rightarrow (-1, 1]$ satisfying the structure condition (\ref{structure condition}).
$r: V\rightarrow (0,+\infty)$ is a Thurston's hyperbolic circle packing metric with $r_i\geq R>0$ for any $i\in V$.
Then there exist positive constants $a_1=a_1(\eta, R)$ and  $a_2=a_2(\eta, R)$
such that $a_1\leq A_{ii}\leq a_2$ for any $i\in V$.
\end{lemma}

\proof
We use $q\sim p$ to denote that there exist  constants $C=C(\eta, R)>0$ and $C'=C'(\eta, R)>0$
such that $C p\leq q \leq C'p$.
Then what we need to prove is equivalent to $A_{ii}\sim 1$.

Note that under the condition $r_i\geq R>0$ for any $i\in V$, we have
\begin{equation}\label{equivalent equation}
C_i\sim e^{r_i}, S_i\sim e^{r_i}, \cosh l_{ij}+1\sim \cosh l_{ij}\sim e^{r_i+r_j},
\end{equation}
where $\cosh l_{ij}\sim e^{r_i+r_j}$ follows from Lemma \ref{equi of lij Thurs hyper CP}.
Combining $\eta\in (-1,1]$, the structure condition (\ref{structure condition}), (\ref{derivative of theta Thurs hyper CP})
and (\ref{equivalent equation}), we have
\begin{equation}\label{equivalent derivative}
\begin{aligned}
\frac{\partial \theta_i^{jk}}{\partial u_j}(\cosh l_{ij}-1)
\sim & \frac{e^{r_i+r_j+r_k}}{A_{ijk}}[(1-\eta_{ij}^2)+\gamma_{ijk}+\gamma_{jik}].
 \end{aligned}
\end{equation}
By the hyperbolic cosine law, we have
\begin{equation}\label{A squre}
\begin{aligned}
A_{ijk}^2
=&\sinh^2 l_{ij}\sinh^2 l_{ik}(1-\cos^2 \theta_i^{jk})\\
=&\sinh^2 l_{ij}\sinh^2 l_{ik}-(\cosh l_{ij}\cosh l_{ik}-\cosh l_{jk})^2\\
=&(\cosh^2 l_{ij}-1)(\cosh^2 l_{ik}-1)-(\cosh l_{ij}\cosh l_{ik}-\cosh l_{jk})^2\\
=&1+2\cosh l_{ij}\cosh l_{ik}\cosh l_{jk}-\cosh^2 l_{ij}-\cosh^2 l_{ik}-\cosh^2 l_{jk}.
 \end{aligned}
\end{equation}
Submitting (\ref{hyper length in r for Th CP}) into (\ref{A squre}) and by lengthy but direct calculations ( refer to the proof of Lemma 2.4 in \cite{Z1} or Lemma 3.1 in \cite{X1}), we have
\begin{equation*}
\begin{aligned}
A_{ijk}^2
=&2S_i^2S_j^2S_k^2(1+\eta_{ij}\eta_{ik}\eta_{jk})+S_i^2S_j^2(1-\eta_{ij}^2)+S_i^2S_k^2(1-\eta_{ik}^2)+S_j^2S_k^2(1-\eta_{jk}^2)\\
&+2C_jC_kS_i^2S_jS_k\gamma_{ijk}+2C_iC_kS_iS_j^2S_k\gamma_{jik}+2C_iC_jS_iS_jS_k^2\gamma_{kij},
 \end{aligned}
\end{equation*}
which implies
\begin{equation}\label{asym of e^{-2(ri+rj+rk)}A^2}
\begin{aligned}
e^{-2(r_i+r_j+r_k)}A_{ijk}^2
\sim [&1+\eta_{ij}\eta_{ik}\eta_{jk}+\gamma_{ijk}+\gamma_{jik}+\gamma_{kij}\\
&+e^{-2r_i}(1-\eta_{jk}^2)+e^{-2r_j}(1-\eta_{ik}^2)+e^{-2r_k}(1-\eta_{ij}^2)]
 \end{aligned}
\end{equation}
by (\ref{equivalent equation}), $\eta\in (-1,1]$ and the structure condition (\ref{structure condition}).
Note that
\begin{equation}\label{eta positive}
\begin{aligned}
1+\eta_{ij}\eta_{ik}\eta_{jk}+\gamma_{ijk}+\gamma_{jik}+\gamma_{kij}=(1+\eta_{ij})(1+\eta_{ik})(1+\eta_{jk})>0
 \end{aligned}
\end{equation}
by $\eta\in (-1,1]$.
Combining (\ref{asym of e^{-2(ri+rj+rk)}A^2}),  (\ref{eta positive}) and $r_i\geq R>0$,
we have $e^{-2(r_i+r_j+r_k)}A_{ijk}^2\sim 1$,
which implies
$$\frac{\partial \theta_i^{jk}}{\partial u_j}(\cosh l_{ij}-1)
\sim  (1-\eta_{ij}^2)+\gamma_{ijk}+\gamma_{jik}$$
by (\ref{equivalent derivative}).
Similarly,
$\frac{\partial \theta_i^{jl}}{\partial u_j}(\cosh l_{ij}-1)
\sim  (1-\eta_{ij}^2)+\gamma_{ijl}+\gamma_{jil}.$
Therefore,
\begin{equation}\label{Aii equivalent}
\begin{aligned}
A_{ii}
\sim & \sum_{j\sim i}[2(1-\eta_{ij}^2)+\gamma_{ijk}+\gamma_{jik}+\gamma_{ijl}+\gamma_{jil}]
 \end{aligned}
\end{equation}
by (\ref{Aii expression}).
By $\eta\in (-1, 1]$ and the structure condition (\ref{structure condition}), we have
$$\sum_{j\sim i}[2(1-\eta_{ij}^2)+\gamma_{ijk}+\gamma_{jik}+\gamma_{ijl}+\gamma_{jil}]\geq 0,$$
where the equality is attained if and only if
$1-\eta_{ij}^2=0$ and $\gamma_{ijk}=\gamma_{jik}=\gamma_{ijl}=\gamma_{jil}=0$
for every vertex $j$ adjacent to $i$.
By $1-\eta_{ij}^2=0$ and $\eta\in (-1, 1]$, we have $\eta_{ij}=1$ for every $j\in V$ adjacent to $i$, which further implies
$\gamma_{jik}=\eta_{ik}+\eta_{ij}\eta_{jk}=1+\eta_{jk}>0$.
This contradicts to $\gamma_{jik}=0$.
Therefore, we have
$\sum_{j\sim i}[2(1-\eta_{ij}^2)+\gamma_{ijk}+\gamma_{jik}+\gamma_{ijl}+\gamma_{jil}]> 0,$
which implies $A_{ii}\sim 1$ by (\ref{Aii equivalent}). This completes the proof.
\qed

As a corollary of  Lemma \ref{estimate of A-ii}, we have the following estimate on the eigenvalues of $L=(\frac{\partial K}{\partial u})$.
\begin{corollary}\label{estimate of eigenvalues of L}
Suppose $(M, \mathcal{T}, \varepsilon, \eta)$ is a weighted triangulated connected closed surface with the weights
$\varepsilon\equiv 1$ and $\eta: E\rightarrow (-1, 1]$ satisfying the structure condition (\ref{structure condition}).
$r: V\rightarrow (0,+\infty)$ is a Thurston's hyperbolic circle packing metric with $r_i\geq R>0$ for any $i\in V$.
Then there exist positive constants $a_3=a_3(\eta, R)$ and  $a_4=a_4(\eta, R)$
such that the eigenvalues $\lambda_1, \cdots, \lambda_N$ of $L=(\frac{\partial K}{\partial u})$
stay in a closed interval  $[a_3, a_4]\subset(0, +\infty)$.
\end{corollary}
\proof
Combining Lemma \ref{lower bound of ch lij-1}, (\ref{theta deriv nonnegative Bij nonpositive})
and Glickenstein-Thomas's formula (\ref{GT's formula}), we have
$$\frac{\partial }{\partial u_i} Area(\triangle ijk)
\geq C_0(\frac{\partial \theta_j^{ik}}{\partial u_i}+\frac{\partial \theta_k^{ij}}{\partial u_i})
=C_0(\frac{\partial \theta_i^{jk}}{\partial u_j}+\frac{\partial \theta_i^{jk}}{\partial u_k})$$
for some positive constant $C_0=C_0(\eta, R)$,
which implies that
\begin{equation}\label{Aii geq Bij}
\begin{aligned}
A_{ii}=\sum_{\triangle ijk\in F}\frac{\partial }{\partial u_i} Area(\triangle ijk) \geq C_0\sum_{j\sim i}(-B_{ij}).
 \end{aligned}
\end{equation}
By (\ref{Aii geq Bij}) and the definition of $B$ in (\ref{defin of B}), the matrix $3A-C_0B$
is diagonal dominant and thus positive definite, which implies $B<\frac{3}{C_0}A$.
Note that $B$ is positive semi-definite by $B_{ij}\leq 0$ for $j\sim i$ in (\ref{theta deriv nonnegative Bij nonpositive}),
we have
$$A\leq L=A+B< (1+\frac{3}{C_0})A.$$
Then the result in the corollary follows by Lemma \ref{estimate of A-ii}.
\qed

Using Corollary \ref{estimate of eigenvalues of L}, we prove the following key lemma on the matrix $L^s$.
\begin{lemma}\label{inequality of L^s_ii}
Suppose $(M, \mathcal{T}, \varepsilon, \eta)$ is a weighted triangulated connected closed surface with the weights
$\varepsilon\equiv 1$ and $\eta: E\rightarrow (-1, 1]$ satisfying the structure condition (\ref{structure condition}).
$r: V\rightarrow (0,+\infty)$ is a Thurston's hyperbolic circle packing metric with $r_i\geq R>0$ for any $i\in V$.
Then for any $s\in \mathbb{R}$, there exists a constant $C=C(s,\eta, R)>0$ such that
\begin{equation*}
\begin{aligned}
\frac{\sum_{j\in V, j\neq i}((L^s)_{ij})^2}{((L^s)_{ii})^2}\leq C \frac{\sum_{j\in V, j\neq i}(L_{ij})^2}{(L_{ii})^2}
 \end{aligned}
\end{equation*}
for all $i\in V$.
\end{lemma}

\proof
By the definition of $L^s=(\frac{\partial K}{\partial u})^s$ and Lemma \ref{positive defi of Ls} (b),
there exists an orthonormal matrix $P\in O(N)$ such that
\begin{equation*}
\begin{aligned}
L^s=P^T\cdot diag\{\lambda_1^s, \cdots, \lambda_n^s\}\cdot P=P^T\cdot \Lambda^s\cdot P,
\end{aligned}
\end{equation*}
where $\lambda_i>0, i=1,\cdots, N,$ are eigenvalues of $L=(\frac{\partial K}{\partial u})$
and $\Lambda=diag\{\lambda_1, \cdots, \lambda_n\}$.
Assume $P=(P_1,\cdots, P_N)$, where $P_1,\cdots, P_N$ are orthonormal column vectors.
Then
$(L^s)_{ij}=P_i^T\Lambda^s P_j,$
which implies
\begin{equation*}
\begin{aligned}
\sum_{j=1}^N((L^s)_{ij})^2
=(L^sL^s)_{ii}
=(P^T\Lambda^{2s}P)_{ii}
=P_i^T\Lambda^{2s}P_i.
\end{aligned}
\end{equation*}
As $P_i$ is a column unit vector,
assume $P_i=(x_1, \cdots, x_N)^T\in \mathbb{R}^N$ with $\sum_{j=1}^Nx_j^2=1$.
Then
\begin{equation*}
\begin{aligned}
\sum_{j=1}^N((L^s)_{ij})^2=\sum_{j=1}^N\lambda_j^{2s}x_j^2,\ \
(L^s)_{ii}=\sum_{j=1}^N\lambda_j^sx_j^2.
\end{aligned}
\end{equation*}
Therefore, we just need to prove that there exists $C=C(s,\eta, R)>0$ such that
\begin{equation}\label{equivalent inequality}
\begin{aligned}
\frac{\sum_{j=1}^N\lambda_j^{2s}x_j^2}{(\sum_{j=1}^N\lambda_j^sx_j^2)^2}-1
\leq C \left(\frac{\sum_{j=1}^N\lambda_j^{2}x_j^2}{(\sum_{j=1}^N\lambda_jx_j^2)^2}-1\right)
 \end{aligned}
\end{equation}
for a unit vector $(x_1, \cdots, x_N)^T\in \mathbb{R}^N$.
Note that
\begin{equation}\label{estimate 1 in proof}
\begin{aligned}
\frac{\sum_{j=1}^N\lambda_j^{2s}x_j^2}{(\sum_{j=1}^N\lambda_j^sx_j^2)^2}-1
=&\frac{(\sum_{j=1}^N\lambda_j^{2s}x_j^2)(\sum_{j=1}^Nx_j^2)-(\sum_{j=1}^N\lambda_j^sx_j^2)^2}
{(\sum_{j=1}^N\lambda_j^sx_j^2)^2}\\
=&\frac{\sum_{ j\neq k}(\lambda_j^{2s}x_j^2x_{k}^2+\lambda_k^{2s}x_k^2x_{j}^2-2\lambda_j^s\lambda_k^sx_j^2x_k^2)}
{2(\sum_{j=1}^N\lambda_j^sx_j^2)^2}\\
=&\frac{\sum_{ j\neq k}(\lambda_j^{s}-\lambda_k^s)^2x_j^2x_{k}^2}
{2(\sum_{j=1}^N\lambda_j^sx_j^2)^2}\\
\leq& \frac{1}{2a_3^{2s}}\sum_{ j\neq k}(\lambda_j^{s}-\lambda_k^s)^2x_j^2x_{k}^2,
 \end{aligned}
\end{equation}
where Corollary \ref{estimate of eigenvalues of L} is used in the last line.
Similarly, we have
\begin{equation}\label{estimate 2 in proof}
\begin{aligned}
\frac{\sum_{j=1}^N\lambda_j^{2}x_j^2}{(\sum_{j=1}^N\lambda_jx_j^2)^2}-1
=\frac{\sum_{ j\neq k}(\lambda_j-\lambda_k)^2x_j^2x_{k}^2}
{2(\sum_{j=1}^N\lambda_jx_j^2)^2}
\geq \frac{1}{2a_4^2}\sum_{ j\neq k}(\lambda_j-\lambda_k)^2x_j^2x_{k}^2.
 \end{aligned}
\end{equation}
By the mean value theorem, there exists $\xi$ between $\lambda_j$ and $\lambda_k$ such that
$\lambda_j^s-\lambda_k^s=s\xi^{s-1}(\lambda_j-\lambda_k)$,
which implies that there exists a constant $C'=C'(s,\eta, R)>0$ such that
\begin{equation}\label{estimate 3 in proof}
\begin{aligned}
|\lambda_j^s-\lambda_k^s|\leq C'|\lambda_j-\lambda_k|
 \end{aligned}
\end{equation}
by Corollary \ref{estimate of eigenvalues of L}.
Then (\ref{equivalent inequality}) is a direct consequence of  (\ref{estimate 1 in proof}), (\ref{estimate 2 in proof}) and
(\ref{estimate 3 in proof}).
\qed

As an application of Lemma \ref{inequality of L^s_ii}, we prove the following comparison between $(L^s)_{ii}$ and
$\sum_{j\in V, j\neq i}|(L^s)_{ij}|$.
\begin{lemma}\label{estimate of L^s}
Suppose $(M, \mathcal{T}, \varepsilon, \eta)$ is a weighted triangulated connected closed surface with the weights
$\varepsilon\equiv 1$ and $\eta: E\rightarrow (-1, 1]$ satisfying the structure condition (\ref{structure condition}).
$r: V\rightarrow (0,+\infty)$ is a Thurston's hyperbolic circle packing metric with $r_i\geq R>0$ for any $i\in V$.
For any $s\in \mathbb{R}$ and $\widetilde{C}>0$, there exists a constant $R_1=R_1(s, \eta, R, \widetilde{C})>0$
such that if $r_i>R_1$, then
$$(L^s)_{ii}\geq \widetilde{C}\sum_{j\in V, j\neq i}|(L^s)_{ij}|.$$
\end{lemma}
\proof
By Lemma \ref{equi of lij Thurs hyper CP},
there exists $C_1=C_1(\eta)>0$ such that
$\cosh l_{ij}\geq C_1(C_iC_j+S_iS_j)=C_1\cosh (r_i+r_j)> C_1\cosh r_i$.
Therefore, for any $C_2>0$ (to be determined), there exists $R_1=R_1(\eta, C_1, C_2)>0$ such that if $r_i>R_1$,
then $\cosh l_{ij}-1\geq C_2$
for every edge $\{ij\}\in E$ adjacent to $i$, which implies
\begin{equation}\label{Aii geq C sum Bij}
\begin{aligned}
A_{ii}\geq C_2\sum_{j\sim i}\left(\frac{\partial \theta_j^{ik}}{\partial u_i}+\frac{\partial \theta_k^{ij}}{\partial u_i}\right)
=C_2\sum_{j\sim i}(-B_{ij})
\end{aligned}
\end{equation}
by (\ref{Aii expression}).
Therefore,
$$L_{ii}=A_{ii}+\sum_{j\sim i}(-B_{ij})\geq (C_2+1)\sum_{j\sim i}(-B_{ij})=(C_2+1)\sum_{j\in V, j\neq i}|L_{ij}|,$$
by (\ref{Aii geq C sum Bij}),
which further implies
\begin{equation}\label{Lii squre geq C sum Lij squre}
\begin{aligned}
L_{ii}^2\geq (C_2+1)^2\sum_{j\in V, j\neq i}|L_{ij}|^2.
\end{aligned}
\end{equation}
Combining (\ref{Lii squre geq C sum Lij squre}) and Lemma \ref{inequality of L^s_ii}, we have
\begin{equation}\label{proof inequality Lii squre}
\begin{aligned}
\frac{\sum_{j\in V, j\neq i}((L^s)_{ij})^2}{((L^s)_{ii})^2}\leq C_3\frac{\sum_{j\in V, j\neq i}(L_{ij})^2}{(L_{ii})^2}\leq \frac{C_3}{(C_2+1)^2},
\end{aligned}
\end{equation}
where $C_3=C_3(s,\eta, R)>0$ is given by Lemma \ref{inequality of L^s_ii}.
The inequality (\ref{proof inequality Lii squre}) implies
\begin{equation*}
\begin{aligned}
(L^s)_{ii}\geq \frac{C_2+1}{\sqrt{C_3}} |(L^s)_{ij}|
\end{aligned}
\end{equation*}
for $j\neq i$ and then
\begin{equation*}
\begin{aligned}
(L^s)_{ii}\geq \frac{C_2+1}{N\sqrt{C_3}}\sum_{j\in V, j\neq i} |(L^s)_{ij}|.
\end{aligned}
\end{equation*}
Set $C_2=N\widetilde{C}\sqrt{C_3}-1>0$.  Then if $r_i\geq R_1=R_1(\eta, C_1, C_2)=R_1(s,\eta, R, \widetilde{C})$,
we have $(L^s)_{ii}\geq \widetilde{C}\sum_{j\in V, j\neq i}|(L^s)_{ij}|$.
\qed

\begin{remark}
In the case of $s=1$ and $\eta\in [0,1]$, the result in Lemma \ref{estimate of L^s}
was proved by Ge-Hua \cite{GH1}.
\end{remark}

We shall prove that the solution $u(t)$ of fractional combinatorial Calabi flow (\ref{FCCF equation})
for Thurston's hyperbolic circle packings
stays in a compact subset of $\mathbb{R}^N_{<0}$.
To prove this, we further need the following result on Thurston's hyperbolic circle packing metrics.
\begin{lemma}\label{theta small for ri large}
Suppose $(M, \mathcal{T}, \varepsilon, \eta)$ is a weighted triangulated connected closed surface with the weights
$\varepsilon\equiv 1$ and $\eta: E\rightarrow (-1, 1]$ satisfying the structure condition (\ref{structure condition}).
$\triangle ijk$ is a triangle in $F$.
For any $\epsilon>0$, there exists a positive constant $R_2$ such that if $r_i>R_2$, then $\theta_i^{jk}<\epsilon$.
\end{lemma}
Lemma \ref{theta small for ri large} has been proved for different cases in almost the same manner.
We will not give another proof for Lemma \ref{theta small for ri large} here.
Readers interested in the proof could refer to \cite{GJ1, GX4,GX5,GHZ, S, X9, Z1}.

\begin{proposition}\label{upper bound of u_i(t) for hyperbolic Thurs CP}
Suppose $(M, \mathcal{T}, \varepsilon, \eta)$ is a weighted triangulated connected closed surface with the weights
$\varepsilon\equiv 1$ and $\eta: E\rightarrow (-1, 1]$ satisfying the structure condition (\ref{structure condition}).
Suppose that there exists  $\overline{u}\in \mathbb{R}^N_{<0}$ such that $K(\overline{u})=\overline{K}$
in the hyperbolic background geometry and
$u(t), t\in [0, T), $ is a solution of the hyperbolic fractional combinatorial Calabi flow (\ref{FCCF equation})
with the maximal existing time $T\leq +\infty$.
Then there exists a constant $c<0$ such that $u_i(t)<c$ for every $i\in V$
along the hyperbolic fractional combinatorial Calabi flow (\ref{FCCF equation}).
\end{proposition}
\proof
By the relationship of $u_i$ and $r_i$ in (\ref{ui ri relationship}), we just need to prove that $r_i(t)$ is uniformly bounded
from above in $(0, +\infty)$ for every $i\in V$ along the hyperbolic fractional combinatorial Calabi flow (\ref{FCCF equation}),
which is equivalent to $r(t)$ is bounded along (\ref{FCCF equation}).

By Lemma \ref{lower bound of u_i(t) for hyperbolic Thurs CP} and Remark \ref{lower bound of r_i geq R},
there exists a constant $R>0$ such that $r_i(t)\geq R>0$ along the fractional combinatorial Calabi flow (\ref{FCCF equation})
for every $i\in V$ under the existence of $\overline{u}$ with $K(\overline{u})=\overline{K}$.
By Lemma \ref{theta small for ri large}, for $\epsilon=\frac{1}{2N}(2\pi-\overline{K}_i)>0$,
there exists a constant $R_2>0$ with $R_2\geq R$ such that if $r_i>R_2$,
then $\theta_{i}^{jk}<\epsilon$ for every triangle $\triangle ijk\in F$ at $i$,
which implies
\begin{equation}\label{proof inequ 1 in main propo}
\begin{aligned}
K_i-\overline{K}_i=2\pi-\sum_{\triangle ijk\in F} \theta_{i}^{jk}-\overline{K}_i> \frac{1}{2}(2\pi-\overline{K}_i).
\end{aligned}
\end{equation}
By Lemma \ref{estimate of L^s}, for $\widetilde{C}=\frac{2\max_{j\in V}[(N+2)\pi+|\overline{K}_j|]}{(2\pi-\overline{K}_i)}>0$,
there exists $R_1>0$ with $R_1\geq R$ such that if $r_i\geq R_1$, then
\begin{equation}\label{proof inequ 2 in main propo}
\begin{aligned}
(L^s)_{ii}\geq \widetilde{C}\sum_{j\in V, j\neq i}|(L^s)_{ij}|
=\frac{2\max_{j\in V}[(N+2)\pi+|\overline{K}_j|]}{(2\pi-\overline{K}_i)}\sum_{j\in V, j\neq i}|(L^s)_{ij}|.
\end{aligned}
\end{equation}
Set $R_3= \max (R_1, R_2, r_i(0)+1)$.  If $r_i\geq R_3$,  then
\begin{equation}\label{dui dt less 0}
\begin{aligned}
L^s(K-\overline{K})_i
=&(L^s)_{ii}(K_i-\overline{K}_i)+\sum_{j\in V, j\neq i}(L^s)_{ij}(K_j-\overline{K}_j)\\
> &\frac{1}{2}(2\pi-\overline{K}_i)(L^s)_{ii}+\sum_{j\in V, j\neq i}(L^s)_{ij}(K_j-\overline{K}_j)\\
\geq & \max_{j\in V}[(N+2)\pi+|\overline{K}_j|]\sum_{j\in V, j\neq i}|(L^s)_{ij}|+\sum_{j\in V, j\neq i}(L^s)_{ij}(K_j-\overline{K}_j)\\
= &\sum_{j\in V, j\neq i}|(L^s)_{ij}|\left(\max_{k\in V}[(N+2)\pi+|\overline{K}_k|]-|K_j-\overline{K}_j|\right)\\
\geq&0,
\end{aligned}
\end{equation}
where (\ref{proof inequ 1 in main propo}) is used in the second line and  (\ref{proof inequ 2 in main propo}) is used in the third line.

Suppose that along the fractional combinatorial Calabi flow (\ref{FCCF equation}), $r(t), t\in [0, T)$, is not bounded.
Then there exists at least one vertex $i_0\in V$ such that $\limsup_{t\rightarrow T^-} r_{i_0}=+\infty$.
Without loss of generality, we can take $i_0=i$.
Therefore, there exists $t_0\in (0, T)$ such that $r_i(t_0)>R_3$ by $\limsup_{t\rightarrow T^-} r_{i}=+\infty$.
Set
$$a=\inf \{t<t_0|r_{i}(s)>R_3, \forall s\in [t, t_0]\}.$$
Then $a\in (0, t_0)$ by $R_3\geq r_i(0)+1$, $r_i(a)=R_3$ and $r_i(t)>R_3$ for all $t\in (a, t_0]$.
Combining this with (\ref{dui dt less 0}),
we have
$$\frac{du_i}{dt}
=-L^s(K-\overline{K})_i<0$$
for $t\in (a, t_0]$ along the hyperbolic fractional combinatorial Calabi flow (\ref{FCCF equation}), which implies $\frac{dr_i}{dt}=\frac{dr_i}{du_i}\frac{du_i}{dt}<0$
and then $r_{i}(t)<r_{i}(a)=R_3$ for $t\in (a, t_0]$.
This contradicts to the fact that $r_i(t_0)>R_3$.
Therefore, $\limsup_{t\rightarrow T^-} r_{i}<+\infty$ for every $i\in V$,
which implies $r_i(t)$ is bounded from above for every $i\in V$.
\qed

As a direct corollary of Lemma \ref{lower bound of u_i(t) for hyperbolic Thurs CP} and Proposition \ref{upper bound of u_i(t) for hyperbolic Thurs CP}, we have the following result,
which proves the longtime existence part of Theorem \ref{main theorem conv of FCCF Thurston's CP} (b).

\begin{corollary}
Suppose $(M, \mathcal{T}, \varepsilon, \eta)$ is a weighted triangulated connected closed surface with the weights
$\varepsilon\equiv 1$ and $\eta: E\rightarrow (-1, 1]$ satisfying the structure condition (\ref{structure condition}).
Suppose there exists  $\overline{u}\in \mathbb{R}^N_{<0}$ such that $K(\overline{u})=\overline{K}$
in the hyperbolic background geometry.
Then the solution  $u(t)$ of the hyperbolic fractional combinatorial Calabi flow (\ref{FCCF equation})
stays in a compact subset of $\mathbb{R}^N_{<0}$ and exists for all time.
\end{corollary}

The following of the proof for Theorem \ref{main theorem conv of FCCF Thurston's CP} (b), i.e. the convergence of the solution
$u(t)$ of  hyperbolic fractional combinatorial Calabi flow (\ref{FCCF equation}) to $\overline{u}$
in the case of Thurston's hyperbolic circle packing metrics,
is almost the same as that for Theorem \ref{main theorem conv of FCCF Thurston's CP} (a). We omit the details here.

\section{Fractional combinatorial Calabi flow for vertex scaling on surfaces}\label{section 4}

Vertex scaling is a special type of discrete conformal structures in Definition \ref{definition of DCS} with $\varepsilon\equiv0$.
In the Euclidean background geometry, $\varepsilon\equiv0$ implies $l_{ij}^2=\eta_{ij}e^{f_i+f_j}$
by Definition \ref{definition of DCS}, which implies that
\begin{equation}\label{Euclidean VS defi}
\begin{aligned}
\widetilde{l}_{ij}=\l_{ij}e^{\frac{u_i+u_j}{2}}
\end{aligned}
\end{equation}
for two different PL metrics $l$ and $\widetilde{l}$ on $(M, \mathcal{T}, \varepsilon, \eta)$
with $u=\widetilde{f}-f$.
(\ref{Euclidean VS defi}) is the definition of Euclidean vertex scaling on surfaces introduced independently by Luo \cite{L1} and
R\v{o}cek-Williams \cite{RW}.
In the hyperbolic background geometry, $\varepsilon\equiv0$ implies $\sinh^2\frac{l_{ij}}{2}=\frac{\eta_{ij}}{2}e^{f_i+f_j}$
by Definition \ref{definition of DCS}, which further implies that
\begin{equation}\label{hyperbolic VS defi}
\begin{aligned}
\sinh \frac{\widetilde{l}_{ij}}{2}=\sinh\frac{l_{ij}}{2}e^{\frac{u_i+u_j}{2}}
\end{aligned}
\end{equation}
for two different PH metrics $l$ and $\widetilde{l}$ on $(M, \mathcal{T}, \varepsilon, \eta)$
with $u=\widetilde{f}-f$.
(\ref{hyperbolic VS defi}) is the definition of hyperbolic vertex scaling on surfaces introduced
by Bobenko-Pinkall-Springborn \cite{BPS}.
By comparing with a fixed piecewise Euclidean or hyperbolic metric $l_0$ generated by vertex scaling on a triangulated surface,
it is easy to check that
(\ref{Euclidean VS defi}) and (\ref{hyperbolic VS defi})  are equivalent to the definition of discrete conformal
structures in Definition \ref{definition of DCS} with $\varepsilon\equiv0$ in the Euclidean and hyperbolic
background geometry respectively.
In the following of this section, we will use (\ref{Euclidean VS defi}) and (\ref{hyperbolic VS defi})
as the definition of Euclidean and hyperbolic vertex scaling respectively.
The function $u: V\rightarrow \mathbb{R}$ in (\ref{Euclidean VS defi}) and (\ref{hyperbolic VS defi})
is a shift of the original discrete conformal factor $f:V\rightarrow \mathbb{R}$ in Definition \ref{definition of DCS}
by a constant function defined on $V$ and also called as a discrete conformal factor.

\subsection{The case of Euclidean vertex scaling}
For the Euclidean vertex scaling in (\ref{Euclidean VS defi}), the triangle inequalities are not always satisfied.
As noted in Section \ref{section 2}, this causes that the fractional combinatorial Calabi flow (\ref{FCCF equation}) for
Euclidean vertex scaling may develop singularities.
To handle the potential singularities along the fractional combinatorial Calabi flow (\ref{FCCF equation}) for
Euclidean vertex scaling, we need to do surgery on the flow (\ref{FCCF equation}) by edge flipping under the Delaunay condition.
To prove rigorously the longtime existence and convergence of the solution of
Euclidean fractional combinatorial Calabi flow with surgery,
we need the discrete conformal theory developed by Gu-Luo-Sun-Wu \cite{GLSW} for Euclidean vertex scaling.
We briefly recall some main results in \cite{GLSW} that we need in this paper. For more details, please refer to
Gu-Luo-Sun-Wu's original work \cite{GLSW}.

Gu-Luo-Sun-Wu \cite{GLSW} did not take the triangulations of $(M, V)$ as intrinsic structures attached to the marked surface $(M, V)$,
analogous to the coordinate charts on smooth manifolds.
The intrinsic structure on $(M, V)$ is the polyhedral metric, which is independent of the triangulations.
The combinatorial curvature $K$ is also an intrinsic quantity for a polyhedral metric on $(M, V)$.
Based on this viewpoint, Gu-Luo-Sun-Wu \cite{GLSW} introduced the following
new definition of discrete conformality for PL metrics on $(M, V)$, which allows the triangulations of $(M, V)$
to be changed by edge flipping under the Delaunay condition.

\begin{definition}[\cite{GLSW} Definition 1.1]\label{GLSW's definition of edcf}
Two PL metrics $d, d'$ on $(M, V)$ are discrete conformal if
there exist sequences of PL metrics $d_1=d$, $\cdots$,  $d_m=d'$
on $(M, V)$ and triangulations $\mathcal{T}_1, \cdots, \mathcal{T}_m$ of
$(M, V)$ satisfying
\begin{description}
  \item[(a)] (Delaunay condition) each $\mathcal{T}_i$ is Delaunay in $d_i$,
  \item[(b)] (Vertex scaling condition) if $\mathcal{T}_i=\mathcal{T}_{i+1}$, there exists a function
  $u:V\rightarrow \mathbb{R}$ so that if $e$ is an edge in $\mathcal{T}_i$ with end points $v$ and $v'$,
  then the lengths $l_{d_{i+1}}(e)$ and  $l_{d_{i}}(e)$ of $e$ in $d_i$ and $d_{i+1}$ are related by
  $l_{d_{i+1}}(e)=l_{d_{i}}(e)e^{\frac{u(v)+u(v')}{2}},$
  \item[(c)] if $\mathcal{T}_i\neq \mathcal{T}_{i+1}$, then $(M, d_i)$ is isometric to $(M, d_{i+1})$
  by an isometry homotopic to identity in $(M, V)$.
\end{description}
\end{definition}
The space of PL metrics discrete conformal to $d$ on $(M, V)$ is called as the discrete conformal class of $d$ and
denoted by $\mathcal{D}(d)$.
Gu-Luo-Sun-Wu \cite{GLSW} proved the following discrete uniformization theorem
for PL metrics with discrete conformality given by Definition \ref{GLSW's definition of edcf}.
\begin{theorem}[\cite{GLSW} Theorem 1.2]\label{Euclidean discrete uniformization}
Suppose $(M, V)$ is a closed connected marked surface and $d$ is a PL metric on $(M, V)$.
Then for any $\overline{K}: V\rightarrow (-\infty, 2\pi)$ with $\sum_{v\in V}\overline{K}(v)=2\pi\chi(M)$,
there exists a PL metric $\overline{d}$, unique up to scaling and isometry homotopic to the identity
 on $(M, V)$, such that $\overline{d}$ is discrete conformal to $d$ and the discrete curvature
 of $\overline{d}$ is $\overline{K}$.
\end{theorem}

To prove Theorem \ref{Euclidean discrete uniformization}, Gu-Luo-Sun-Wu \cite{GLSW}
used the decorated Teichim\"{u}ller space theory established by Penner \cite{Penner}.
Denote the decorated Teichim\"{u}ller space of all equivalence class of decorated hyperbolic metrics on $M-V$ by $T_D(M-V)$
and the Teichim\"{u}ller space of all PL metrics on $(M, V)$ by $T_{PL}(M, V)$.
Gu-Luo-Sun-Wu \cite{GLSW} established the following correspondence between $T_D(M-V)$ and $T_{PL}(M, V)$.

\begin{theorem}[\cite{GLSW} Theorem 4.5, Corollary 4.7]\label{main result of GLSW}
There exists a  $C^1$-diffeomorphism $\mathbf{A}: T_{PL}(M, V)\rightarrow T_D(M-V)$
such that $\mathbf{A}|_{\mathcal{D}(d)}: \mathcal{D}(d)\rightarrow \{p\}\times \mathbb{R}^V_{>0}$
is a $C^1$-diffeomorphism.
\end{theorem}
The construction of the map $\mathbf{A}$ is rather technical and we will not give the details of the construction here.
What we need to use is the following property of the map $\mathbf{A}|_{\mathcal{D}(d)}: \mathcal{D}(d)\rightarrow \{p\}\times \mathbb{R}^V_{>0}$.

\begin{theorem}[\cite{GLSW} Proposition 5.2, \cite{L1} Theorem 1.2, Theorem 2.1, Corollary 2.3]\label{energy function Euclidean}
Set $u_i=\ln w_i$ for $w=(w_1, w_2, \cdots, w_n)\in \mathbb{R}^n_{>0}$
and define
\begin{equation}\label{F extension of K}
\begin{aligned}
\mathbf{F}:\mathbb{R}^n&\rightarrow (-\infty, 2\pi)^n\\
u&\mapsto K_{\mathbf{A}^{-1}(p, w(-u))}.
\end{aligned}
\end{equation}
Then
\begin{description}
  \item[(1)] for any $k\in \mathbb{\mathbb{R}}$, $\mathbf{F}(v+k(1, 1, \cdots, 1))=\mathbf{F}(v)$.
  \item[(2)] there exists a $C^2$-smooth convex function $W=\int \sum_{i=1}^n \mathbf{F}_i(u)du_i: \mathbb{R}^n\rightarrow \mathbb{R}$
so that its gradient $\nabla W$ is $\mathbf{F}$ and the restriction $W: \{u\in \mathbb{R}^n|\sum_{i=1}^nu_i=0\}\rightarrow \mathbb{R}$
is strictly convex.
\end{description}
\end{theorem}
Suppose $\mathcal{T}$ is a triangulation of $(M, V)$,
we use $\Omega^{\mathcal{T}}_{D}(d')$ to denote the admissible spaces
of discrete conformal factors $u$ such that $\mathcal{T}$ is Delaunay for $d'\in \mathcal{D}(d)$.
As pointed out in \cite{GLSW},
$\mathbb{R}^n=\cup_{\mathcal{T}} \Omega^{\mathcal{T}}_{D}(d')$ is an finite analytic cell decomposition of $\mathbb{R}^n$.
This further implies that  $\mathbf{F}$ in (\ref{F extension of K}) is a $C^1$ extension of $K|_{\Omega^{\mathcal{T}}_{D}(d')}$
for $d'\in \mathcal{D}(d)$.
Furthermore, by Theorem \ref{energy function Euclidean},
$(\frac{\partial \mathbf{F}}{\partial u})$ is a positive semi-definite matrix with kernel  $\{t\mathbf{1}|t\in \mathbb{R}\}$.
Therefore, we can define the following fractional
discrete Lapalce operator $\widetilde{\Delta}^s$ for discrete conformal factor $u\in \mathbb{R}^N$.
\begin{definition}\label{defi of s order Eucl frac Lap}
Suppose $(M, V)$ is a marked surface with a PL metric $d$ and $s\in \mathbb{R}$ is a constant.
The $2s$-order Euclidean fractional discrete Laplace operator $\widetilde{\Delta}^s$ for $u\in \mathbb{R}^n=\cup_{\mathcal{T}} \Omega^{\mathcal{T}}_{D}(d')$
is defined to be
\begin{equation}
\begin{aligned}
\widetilde{\Delta}^s=-\left(\frac{\partial \mathbf{F}}{\partial u}\right)^s.
\end{aligned}
\end{equation}
\end{definition}

\begin{remark}
If $s=1$, the Euclidean fractional discrete Laplace operator $\widetilde{\Delta}^s$ in Definition \ref{defi of s order Eucl frac Lap}
 is the discrete Laplace operator introduced for
discrete conformal factors $u\in \mathbb{R}^N$ in \cite{GLSW, ZX}.
By Definition \ref{defi of s order Eucl frac Lap},
the fractional discrete Laplace operator
$\widetilde{\Delta}^s=-(\frac{\partial \mathbf{F}}{\partial u})^s$ defined on
$\mathbb{R}^n=\cup_{\mathcal{T}} \Omega^{\mathcal{T}}_{D}(d')$
is a continuous extension of the fractional discrete Laplace operator $\Delta^s=-(\frac{\partial K}{\partial u})^s$
defined on $\Omega^{\mathcal{T}}_{D}(d')$.
Similar to the case for Euclidean discrete Laplace operator  in \cite{BS, GLSW, ZX},
for a PL metric $d$ on $(M, V)$,
the Euclidean fractional discrete Laplace operator $\widetilde{\Delta}^s$ is
an intrinsic operator in the sense that it is
independent of the Delaunay triangulations of $(M, V)$ for $d$.
\end{remark}

Using the fractional discrete Laplace operator $\widetilde{\Delta}^s$ in Definition \ref{defi of s order Eucl frac Lap},
the Euclidean fractional combinatorial Calabi flow (\ref{FCCF equation}) with surgery could be written as
\begin{equation}\label{FFCF with surgery}
\begin{aligned}
\frac{du_i}{dt}=\widetilde{\Delta}^s(\mathbf{F}-\overline{K}).
\end{aligned}
\end{equation}
As the right hand side of (\ref{FFCF with surgery}) is a $C^1$ function of $u\in \mathbb{R}^N$, the local existence for
the solution of Euclidean fractional combinatorial Calabi flow with surgery (\ref{FFCF with surgery}) follows from
the standard theory in ordinary differential equations.

\textbf{Proof of Theorem \ref{main theorem conv of FCCF vertex scaling} (a): }
By Theorem \ref{Euclidean discrete uniformization}, Theorem \ref{main result of GLSW}
and Theorem \ref{energy function Euclidean}, for $\overline{K}: V\rightarrow (-\infty, 2\pi)$
with $\sum_{i=1}^N \overline{K}_i=2\pi\chi(M)$, there exists $\overline{u}\in \mathbb{R}^N$
up to a shift of a vector $t\mathbf{1}$, $t\in \mathbb{R}$, such that $\mathbf{F}(\overline{u})=\overline{K}$.
Without loss of generality, we assume that $\sum_{i=1}^N \overline{u}_i=\sum_{i=1}^N u_i(0)$,
where $u(0)$ is the initial value of (\ref{FFCF with surgery}).

By Theorem \ref{energy function Euclidean}, the kernel of $\widetilde{\Delta}$ is $\{t\mathbf{1}|t\in \mathbb{R}\}$,
which implies the kernel of $\widetilde{\Delta}^s$ is $\{t\mathbf{1}|t\in \mathbb{R}\}$.
Then we have
\begin{equation*}
\begin{aligned}
\frac{d}{dt}\left(\sum_{i=1}^Nu_i(t)\right)=\sum_{i=1}^N\widetilde{\Delta}^s(\mathbf{F}-\overline{K})_i
=\mathbf{1}^T\widetilde{\Delta}^s(\mathbf{F}-\overline{K})=0
\end{aligned}
\end{equation*}
along (\ref{FFCF with surgery}) by $Ker(\widetilde{\Delta}^s)=\{t\mathbf{1}|t\in \mathbb{R}\}$,
which implies that $\sum_{i=1}^Nu_i(t)$ is invariant along (\ref{FFCF with surgery}).
Without loss of generality, assume that $u(0)\in \Sigma_0=\{u\in \mathbb{R}^N| \sum_{i=1}^N u_i=0\}$.
By Theorem \ref{main result of GLSW}, we can define the following energy function
$$\overline{W}(u)=\int_{\overline{u}}^u\sum_{i=1}^N(\mathbf{F}_i-\overline{K}_i)du_i,$$
which is a $C^2$ smooth convex function defined on $\mathbb{R}^N$ with $\nabla \overline{W}=\mathbf{F}-\overline{K}$.
Furthermore, by $\mathbf{F}(\overline{u})=\overline{K}$ and Theorem \ref{energy function Euclidean},
we have $\nabla_u \overline{W}(\overline{u})=0$ and $\lim_{u\in \Sigma_0, u\rightarrow \infty}\overline{W}(u)=+\infty$.

Along the fractional combinatorial Calabi flow with surgery (\ref{FFCF with surgery}),
\begin{equation}\label{dw dt leq 0 FCCF surgery}
\begin{aligned}
\frac{d \overline{W}(u(t))}{dt}
=\sum_{i=1}^N \nabla_{u_i}\overline{W}\cdot \frac{du_i}{dt}
=-(\mathbf{F}-\overline{K})^T\left(\frac{\partial \mathbf{F}}{\partial u}\right)^s(\mathbf{F}-\overline{K})\leq 0
\end{aligned}
\end{equation}
by the property that
$(\frac{\partial \mathbf{F}}{\partial u})^s$ is positive semi-definite with kernel  $\{t\mathbf{1}|t\in \mathbb{R}\}$,
which implies that $\overline{W}(u(t))\leq \overline{W}(u(0))$.
Combining the fact that $u(t)\in \Sigma_0$ along (\ref{FFCF with surgery}) and
$\lim_{u\in \Sigma_0, u\rightarrow \infty}\overline{W}(u)=+\infty$,
this further implies that the solution $u(t)$ of the Euclidean fractional combinatorial Calabi flow with surgery (\ref{FFCF with surgery})
stays in a compact subset of $\Sigma_0$.
Therefore, the solution $u(t)$ of  (\ref{FFCF with surgery}) exists for all time.
The following of the proof is the same as that for Theorem \ref{main theorem conv of FCCF Thurston's CP} (a).
We omit the details here.
\qed

\subsection{The case of hyperbolic vertex scaling}
For hyperbolic vertex scaling, we also need to do surgery by edge flipping under the hyperbolic Delaunay condition
to handle the potential singularities along the fractional combinatorial Calabi flow (\ref{FCCF equation}).
To prove the longtime existence and global convergence for the solution of fractional combinatorial Calabi flow with surgery
in the case of hyperbolic vertex scaling, we need the discrete conformal theory established by
Gu-Guo-Luo-Sun-Wu in \cite{GGLSW}.
Parallelling to the case of Euclidean vertex scaling,
we only sketch the main results obtained by Gu-Guo-Luo-Sun-Wu in \cite{GGLSW} for PH metrics
that we need to use in this paper.
For more details, please refer to  Gu-Guo-Luo-Sun-Wu's original work \cite{GGLSW}.

\begin{definition}[\cite{GGLSW}, Definition 1]\label{hyperbolic discrete conformal}
Two PH metrics $d$, $d'$ on a closed marked surface
$(M, V)$ are discrete conformal if there exists sequences of PH
metrics $d_1=d$, $d_2, \cdots, d_m=d'$ on $(M, V)$ and triangulations $\mathcal{T}_1, \cdots, \mathcal{T}_m$
of $(M, V)$ satisfying
\begin{description}
  \item[(a)] (Delaunay condition)  each $\mathcal{T}_i$ is Delaunay in $d_i$,
  \item[(b)] (Vertex scaling condition) if $\mathcal{T}_i=\mathcal{T}_{i+1}$, there exists a function $u:V\rightarrow \mathbb{R}$,
  called a conformal factor, so that if $e$ is an edge in $\mathcal{T}_i$ with end points $v$ and $v'$,
  then the lengths $x_{d_i}(e)$ and $x_{d_{i+1}}(e)$ of $e$ in metrics $d_i$ and $d_{i+1}$ are related by
  $\sinh \frac{x_{d_{i+1}}(e)}{2}=e^{\frac{u(v)+u(v')}{2}}\sinh \frac{x_{d_{i}}(e)}{2},$
  \item[(c)] if $\mathcal{T}_i\neq \mathcal{T}_{i+1}$, then $(M, d_i)$ is isometric to $(M, d_{i+1})$ by an
  isometry homotopic to the identity in $(M, V)$.
\end{description}
\end{definition}

The space of PH metrics on $(M, V)$ discrete conformal to $d$ is called as the conformal class of $d$ and
denoted by $\mathcal{D}(d)$.
Gu-Guo-Luo-Sun-Wu \cite{GGLSW} proved the following discrete uniformization theorem
for PH metrics with discrete conformality given by Definition \ref{hyperbolic discrete conformal}.
\begin{theorem}[\cite{GGLSW}, Theorem 3]\label{hyperbolic uniformization}
Suppose $(M, V)$ is a closed connected surface with marked points and
$d$ is a PH metric on $(M, V)$.
Then for any $\overline{K}: V\rightarrow (-\infty, 2\pi)$ with
$\sum_{v\in V}\overline{K}(v)>2\pi\chi(M)$, there exists a unique
PH metric $d'$ on $(M, V)$ so that $d'$ is discrete conformal
to $d$ and the discrete curvature of $d'$ is $\overline{K}$.
\end{theorem}

Following the Euclidean case in \cite{GLSW}, Gu-Guo-Luo-Sun-Wu \cite{GGLSW}
denote the decorated Teichim\"{u}ller space of all equivalence class of decorated hyperbolic metrics on $M-V$ by $T_D(M-V)$
and the Teichim\"{u}ller space of all PH metrics on $(M, V)$ by $T_{hp}(M, V)$.
Based on Penner's decorated Teichim\"{u}ller space theory \cite{Penner},
Gu-Guo-Luo-Sun-Wu \cite{GGLSW} established the following correspondence between $T_D(M-V)$ and $T_{hp}(M, V)$.

\begin{theorem}[\cite{GGLSW} Theorem 22, Corollary 24]\label{main result of GGLSW}
There exists a  $C^1$-diffeomorphism $\mathbf{A}: T_{hp}(M, V)\rightarrow T_D(M-V)$
such that $\mathbf{A}|_{\mathcal{D}(d)}: \mathcal{D}(d)\rightarrow \{p\}\times \mathbb{R}^N_{>0}$
is a $C^1$-diffeomorphism.
\end{theorem}
Parallelling to the Euclidean case, Gu-Guo-Luo-Sun-Wu \cite{GGLSW}
further prove the following property of the map $\mathbf{A}|_{\mathcal{D}(d)}: \mathcal{D}(d)\rightarrow \{p\}\times \mathbb{R}^N_{>0}$ in the hyperbolic background geometry.

\begin{theorem}[\cite{GGLSW} Section 4.1, \cite{BPS} Proposition 6.1.5]\label{energy function hyperbolic}
Set $u_i=\ln w_i$ for $w=(w_1, w_2, \cdots, w_N)\in \mathbb{R}^N_{>0}$
and define
\begin{equation}\label{F extension of K hyper}
\begin{aligned}
\mathbf{F}:\mathbb{R}^n&\rightarrow (-\infty, 2\pi)^N\\
u&\mapsto K_{\mathbf{A}^{-1}(p, w(-u))}.
\end{aligned}
\end{equation}
Then there exists a $C^2$-smooth strictly convex function $W=\int \sum_{i=1}^N \mathbf{F}_i(u)du_i: \mathbb{R}^n\rightarrow \mathbb{R}$ with gradient $\nabla W=\mathbf{F}$.
\end{theorem}

Similar to the Euclidean case, suppose $\mathcal{T}$ is a triangulation of $(M, V)$,
we use $\Omega^{\mathcal{T}}_{D}(d')$ to denote the admissible spaces
of hyperbolic discrete conformal factors $u$ such that $\mathcal{T}$ is Delaunay for $d'\in \mathcal{D}(d)$.
Then  $\mathbb{R}^N=\cup_{\mathcal{T}} \Omega^{\mathcal{T}}_{D}(d')$
is an finite analytic cell decomposition of $\mathbb{R}^n$ \cite{GGLSW}.
By Theorem \ref{main result of GGLSW} and Theorem \ref{energy function hyperbolic},
$\mathbf{F}$ is a $C^1$ smooth function defined on $\mathbb{R}^N=\cup_{\mathcal{T}} \Omega^{\mathcal{T}}_{D}(d')$
and $(\frac{\partial \mathbf{F}}{\partial u})$ is a continuous strictly positive definite matrix defined for  $u\in \mathbb{R}^N$.
Specially, this implies that $(\frac{\partial \mathbf{F}}{\partial u})$ is independent of the Delaunay triangulation of $(M, V)$
for a PH metric $d'\in \mathcal{D}(d)$.
Following the Euclidean case, we can define the following hyperbolic fractional discrete Laplace operator.

\begin{definition}\label{defi of s order hyper frac Lap}
Suppose $(M, V)$ is a marked surface with  a PH metric $d$ and $s\in \mathbb{R}$ is a constant.
The $2s$-order hyperbolic fractional discrete Laplace operator $\widetilde{\Delta}^s$ for $u\in \mathbb{R}^n=\cup_{\mathcal{T}} \Omega^{\mathcal{T}}_{D}(d')$
is defined to be
\begin{equation}
\begin{aligned}
\widetilde{\Delta}^s=-\left(\frac{\partial \mathbf{F}}{\partial u}\right)^s.
\end{aligned}
\end{equation}
\end{definition}

\begin{remark}
If $s=1$, the hyperbolic fractional discrete Laplace operator $\widetilde{\Delta}^s$ in Definition \ref{defi of s order hyper frac Lap}
 is the hyperbolic discrete Laplace operator introduced for
discrete conformal factors $u\in \mathbb{R}^N$ in \cite{GGLSW, ZX}.
By Definition \ref{defi of s order hyper frac Lap},
the fractional discrete Laplace operator
$\widetilde{\Delta}^s=-(\frac{\partial \mathbf{F}}{\partial u})^s$ defined on
$\mathbb{R}^n=\cup_{\mathcal{T}} \Omega^{\mathcal{T}}_{D}(d')$
is a continuous extension of the fractional discrete Laplace operator $\Delta^s=-(\frac{\partial K}{\partial u})^s$
defined on $\Omega^{\mathcal{T}}_{D}(d')$.
Similar to the hyperbolic discrete Laplace operator in \cite{GGLSW, ZX},
for a PH metric $d$ on $(M, V)$,
the hyperbolic fractional discrete Laplace operator $\widetilde{\Delta}^s$ is
an intrinsic operator for $d$ in the sense that it is
independent of the Delaunay triangulations of $(M, V)$ for $d$.
\end{remark}

Using the hyperbolic fractional discrete Laplace operator $\widetilde{\Delta}^s$ in Definition \ref{defi of s order hyper frac Lap},
the fractional combinatorial Calabi flow (\ref{FCCF equation}) with surgery for hyperbolic vertex scaling could be written as
\begin{equation}\label{FFCF with surgery hyperbolic}
\begin{aligned}
\frac{du_i}{dt}=\widetilde{\Delta}^s(\mathbf{F}-\overline{K}).
\end{aligned}
\end{equation}
As the right hand side of (\ref{FFCF with surgery hyperbolic}) is a continuous function of $u\in \mathbb{R}^N$,
the local existence for the solution of hyperbolic fractional combinatorial Calabi flow
with surgery (\ref{FFCF with surgery hyperbolic}) follows from
the standard theory of ordinary differential equations.

\textbf{Proof of Theorem \ref{main theorem conv of FCCF vertex scaling} (b): }
By Theorem \ref{hyperbolic uniformization}, Theorem \ref{main result of GGLSW} and Theorem \ref{energy function hyperbolic},
for $\overline{K}: V\rightarrow (-\infty, 2\pi)$
with $\sum_{i=1}^N \overline{K}_i>2\pi\chi(M)$, there exists a unique $\overline{u}\in \mathbb{R}^N$  such that $\mathbf{F}(\overline{u})=\overline{K}$.

By Theorem \ref{energy function hyperbolic}, we can define the following energy function
$$\overline{W}(u)=\int_{\overline{u}}^u\sum_{i=1}^N(\mathbf{F}_i-\overline{K}_i)du_i,$$
which is a $C^2$ smooth strictly convex function defined on $\mathbb{R}^N$ with $\nabla \overline{W}=\mathbf{F}-\overline{K}$.
Furthermore, by $\mathbf{F}(\overline{u})=\overline{K}$ and Theorem \ref{energy function hyperbolic},
we have $\nabla_u \overline{W}(\overline{u})=0$ and $\lim_{u\rightarrow \infty}\overline{W}(u)=+\infty$.

Along the hyperbolic fractional combinatorial Calabi flow with surgery (\ref{FFCF with surgery hyperbolic}),
we have
\begin{equation}\label{dw dt leq 0 FCCF surgery hyper}
\begin{aligned}
\frac{d \overline{W}(u(t))}{dt}
=\sum_{i=1}^N \nabla_{u_i}\overline{W}\cdot \frac{du_i}{dt}
=-(\mathbf{F}-\overline{K})^T\left(\frac{\partial \mathbf{F}}{\partial u}\right)^s(\mathbf{F}-\overline{K})\leq 0
\end{aligned}
\end{equation}
by the property that
$(\frac{\partial \mathbf{F}}{\partial u})^s$ is a strictly positive definite matrix,
which implies that $\overline{W}(u(t))\leq \overline{W}(u(0))$ along (\ref{FFCF with surgery hyperbolic}).
Combining with
$\lim_{u\rightarrow \infty}\overline{W}(u)=+\infty$,
this further implies that the solution $u(t)$ of hyperbolic
fractional combinatorial Calabi flow with surgery (\ref{FFCF with surgery hyperbolic})
stays in a compact subset of $\mathbb{R}^N$.
Therefore, the solution $u(t)$ of  (\ref{FFCF with surgery hyperbolic}) exists for all time.
The following of the proof is the same as that for Theorem \ref{main theorem conv of FCCF Thurston's CP} (a).
We omit the details here.
\qed

(Tianqi Wu) Center of Mathematical Sciences and Applications, Harvard University, Cambridge, MA 02138

E-mail: tianqi@cmsa.fas.harvard.edu\\[2pt]

(Xu Xu) School of Mathematics and Statistics, Wuhan University, Wuhan 430072, P.R. China

E-mail: xuxu2@whu.edu.cn\\[2pt]

\end{document}